\documentclass[onecolumn]{article}
\usepackage[preprint]{log_2025}

\makeatletter
\renewcommand{\@noticestring}{}
\makeatother

\usepackage[ruled]{algorithm2e}
\usepackage{graphicx}  
\usepackage{amsfonts}
\usepackage{amsmath} 
\usepackage{booktabs}
\usepackage{siunitx}
\usepackage[table, svgnames, dvipsnames]{xcolor}
\usepackage{float}
\usepackage{accents}
\usepackage{tabularx}
\usepackage{array}
\usepackage{comment}
\usepackage{mathtools}
\usepackage{wrapfig}
\usepackage{picinpar}
\usepackage{pifont}
\usepackage{afterpage}
\usepackage[numbers]{natbib} 
\definecolor{cyan}{HTML}{00AEEF}
\definecolor{royalblue}{HTML}{00B4CE}
\definecolor{wildstrawberry}{HTML}{EE2967}
\usepackage{tikz}
\usepackage{multirow}

\DeclareMathOperator{\EX}{\mathbb{E}}
\usepackage{abstract}

\newcommand{\textBF}[1]{%
    \pdfliteral direct {2 Tr 0.3 w} 
     #1%
    \pdfliteral direct {0 Tr 0 w}%
}

\definecolor{light-gray}{gray}{0.95}
\definecolor{Cerulean}{HTML}{00A2E3}
\definecolor{pinegreen}{HTML}{009B55}
\colorlet{LightCerulean}{Cerulean!30}

\colorlet{Sand}{SpringGreen!45}

\definecolor{SeaGreen}{HTML}{3FBC9D}
\colorlet{SeaGreen}{SpringGreen!90}

\usepackage{tocloft}

\def\Xint#1{\mathchoice
{\XXint\displaystyle\textstyle{#1}}%
{\XXint\textstyle\scriptstyle{#1}}%
{\XXint\scriptstyle\scriptscriptstyle{#1}}%
{\XXint\scriptscriptstyle\scriptscriptstyle{#1}}%
\!\int}
\def\XXint#1#2#3{{\setbox0=\hbox{$#1{#2#3}{\int}$ }
\vcenter{\hbox{$#2#3$ }}\kern-.6\wd0}}

\def\dashint{\Xint-}

\newcolumntype{P}[1]{>{\centering\arraybackslash}p{#1}}

\newcommand*{\belowrulesepcolor}[1]{%
  \noalign{%
    \kern-\belowrulesep 
    \begingroup 
      \color{#1}%
      \hrule height\belowrulesep 
    \endgroup 
  }%
} 
\newcommand*{\aboverulesepcolor}[1]{%
  \noalign{%
    \begingroup 
      \color{#1}%
      \hrule height\aboverulesep 
    \endgroup 
    \kern-\aboverulesep 
  }%
}

\usepackage{fancyhdr}

\makeatletter

\title{Geometric flow regularization in latent spaces for smooth dynamics with the efficient variations of curvature}
\author{Andrew Gracyk\footnotemark[1]}
\date{}

\begin{document}

\vspace{-20mm}
\begin{@twocolumnfalse}
\maketitle

\begin{abstract}
\vspace{0mm}
We design strategies in nonlinear geometric analysis to temper the effects of adversarial learning for sufficiently smooth data of numerical method-type dynamics in encoder-decoder methods, variational and deterministic, through the use of geometric flow regularization. We augment latent spaces with geometric flows to control structure, relying on adaptations of curvature and Ricci flow. All of our flows are solved using physics-informed learning. Traditional geometric meaning is traded for computing ability, but we maintain key geometric invariants, the primary of which are maintained, intrinsically-low structure, nontriviality due to sufficient lower bounds on curvature, distortion of volume element, that develop quality in the inference stage. We instill representations that are canonical, smooth, curvature-aware, geodesic-aware, and non-topologically void or sparse. The primary bottleneck of a Ricci curvature flow is that Ricci curvature is high order, thus expensive to compute, so we will attempt to overcome this with properly justified proxies. Our primary contributions are fourfold. We develop a loss based on Gaussian curvature using closed path circulation integration for surfaces, bypassing automatic differentiation of the Christoffel symbols through use of Stokes' theorem. We invent a new parametric flow valid under a Taylor expansion derived from the Gauss equation. We develop two strategies based on time differentiation of functionals, one with a special case of scalar curvature for conformally-changed metrics, and another with harmonic maps, their energy, and induced metrics. Our losses are diminished analytically and mostly heuristic but maintain overall integral latent structure. We showcase that curvature flows and the formulation of geometric structure in intermediary encoded settings enhance learning and overall zero-shot and adversarial fidelity.
\vspace{4mm}
\end{abstract}
\vspace{2mm}
\end{@twocolumnfalse}

\footnotetext[1]{\normalsize Partially supported by NSF Grant 1922758 hosted by the National Center for Supercomputing Applications. Purdue University, Department of Mathematics, West Lafayette, IN 47907, United States, \texttt{agracyk@purdue.edu}.}

\tableofcontents

\section{Introduction}
\label{sec:introduction}

\begin{wrapfigure}{r}{0.55\textwidth}
  \centering
  \vspace{-5mm}
  \includegraphics[width=\linewidth]{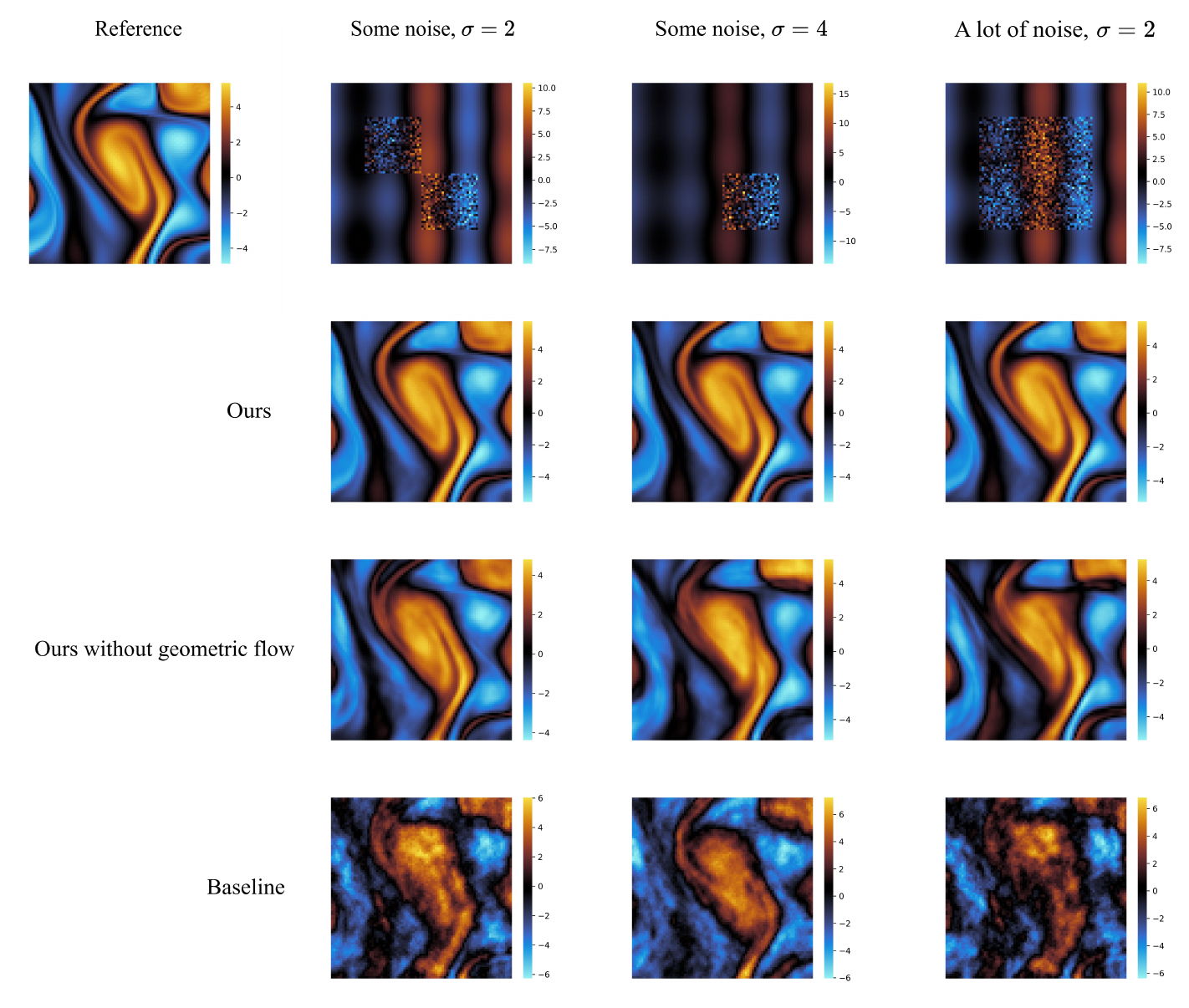}
  \caption{We illustrate scenarios of noised data in our Navier-Stokes experiment.}
  \label{fig:navier_stokes_ood}
  \vspace{-5mm}
\end{wrapfigure}

We will be studying the auxiliary qualities of PDE learning via mitigating the adversarial effects of irregular data, such as that which is out-of-distribution and noised, by availing the development of geometric structure in the latent stages in learning tasks. Our methods will emphasize computational efficiency in the training stage for scaling in moderate intrinsic dimensions primarily, but we will develop techniques that permit high-dimensional representations as well. By controlling the geometries in these intermediary stages, we will study how the development of these structures helps induce accommodating representations and their interconnections with zero-shot generalization and adversarial alleviation. In general, specific geometric qualities in variational and deterministic encoder-decoder settings have been demonstrated to help robustness \cite{lopez2021variationalautoencoderslearningnonlinear} \cite{lopez2025gdvaesgeometricdynamicvariational} \cite{lee2023explicitcurvatureregularizationdeep} \cite{khan2024adversarialrobustnessvaeslens} \cite{camuto2022variationalautoencodersharmonicperspective} \cite{shukla2019geometrydeepgenerativemodels}, inference capacity \cite{nickel2017poincareembeddingslearninghierarchical}, sampling quality \cite{chadebec2020geometryawarehamiltonianvariationalautoencoder} \cite{chadebec2022geometricperspectivevariationalautoencoders} \cite{arvanitidis2021latentspaceodditycurvature} \cite{nagano2019wrappednormaldistributionhyperbolic}, online performance \cite{chen2020learningflatlatentmanifolds}, identifiability \cite{lopez2025gdvaesgeometricdynamicvariational} \cite{arvanitidis2020geometricallyenrichedlatentspaces}, clustering understanding \cite{nagano2019wrappednormaldistributionhyperbolic}, alignment \cite{pmlr-v243-lahner24a}, and anomaly detection \cite{kato2020ratedistortionoptimizationguidedautoencoder}. Our methods are consistent with general strategies of manifold regularization \cite{belkin_manifold}. In particular, our techniques are structured so that initial data is passed into local coordinates before it is mapped to the latent manifold regularizer with an intrinsic curvature flow  \cite{toppingricciflow} \cite{perelman2002entropyformularicciflow} \cite{ricciflowsheridan} \cite{calegariricciflow}.

\vspace{2mm}

\begin{figure}[b]
  \centering
  \includegraphics[width=0.85\linewidth]{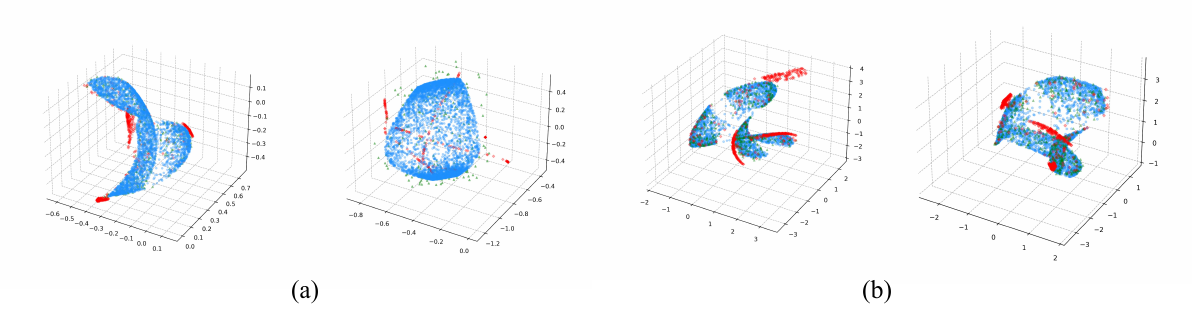}
  \caption{We illustrate the manifold latent space in the extended case where all $u,\mathcal{E},\mathcal{D}$ are active. (a) corresponds to a curvature-regularized method, and (b) corresponds to a latent space that forms naturally. As we can see, ours are more canonical with respect to a curvature-uniform metric. Red corresponds to out-of-distribution data, and green is noised data. These results are consistent with our Gaussian curvature experiment. We remark in much of our experiments, we initialize spherically, thus our flows have consistency with a somewhat curvature-uniform geometry.}
  \label{fig:manifold_figures}
\end{figure}

Our methods are analytically heuristic, and this is intentional since it facilitates efficiency if we diminish theoretical effect. Our proposed methods control structure in ways that (1) prevent contortion; (2) are differentiable, i.e. bypass singularity; (3) have a sense of curvature-uniformization; (4) prevent lapses in point cloud density and void topologies; (5) non-degenerate in measure of the continuously-interpolated manifold; (6) regulate how geodesics exist along the manifold. We refer to \ref{fig:manifold_figures} for an illustration. It is not imperative our methods possesses rigorous analytic qualities that hold up to theoretical expectation, but rather simply that the geometries maintain overall sound structure. It is well known in literature that empirical outcomes are influenced by the manner of representation. \cite{sun2025geometryawaregenerativeautoencoderswarped} \cite{shao2017riemanniangeometrydeepgenerative} \cite{arvanitidis2020geometricallyenrichedlatentspaces} discuss how representational geometries affect learning, and \cite{Picot_2023} \cite{Cosmo_2020} \cite{Yan2024EnhanceAR} \cite{arvanitidis2020geometricallyenrichedlatentspaces} \cite{10.1145/3534678.3539117} \cite{hoffman2019robustlearningjacobianregularization} discuss how certain geometries have regularizing properties via distortion, geodesics, curvature, and void regions. References that discuss how curvature specifically affects adversarial learning or can be used in regularization are \cite{khan2024adversarialrobustnessvaeslens} \cite{moosavidezfooli2018robustnesscurvatureregularizationvice} \cite{Tron_2024} \cite{kim2024gammavaecurvatureregularizedvariational}.

\vspace{2mm}

\section{Methods}
\label{sec:methods}

We will be learning geometric flows of a generalized form $
\partial_t g + \mathcal{F}_{K}[g] = 0$. Ultimately, we seek a series of functions such that $
\tilde{\phi}_t \approx  ( \mathcal{D}\circ \mathcal{E} \circ u) ( \tilde{\phi}_0, t)$.
Consequently, we begin by parameterizing three neural networks $u : \Phi  \times \Theta_u \rightarrow \Sigma \subseteq \mathbb{R}^d, \mathcal{E} : \Sigma \times [0,T] \times \Theta_{\mathcal{E}} \rightarrow \mathcal{M} \subset \mathbb{R}^D, \mathcal{D} : \mathbb{R}^D \times \Theta_{\mathcal{D}} \rightarrow \Phi'$, which act as (1) a manifold parameterizer mapping to local coordinates $\Sigma$; (2) an encoder mapping to the manifold immersed in $\mathbb{R}^D$; (3) a decoder mapping to the PDE solution respectively. We will assume $\Sigma$ is compact throughout the rest of this work. Mapping PDE data directly onto a manifold without local coordinates presents challenges. Generally, a metric as well as tangent spaces are intrinsic, but we need to work with induced metrics. Thus, by parameterizing the manifold, we can perform the necessary calculations. A Ricci-type flow is placed upon the intermediary $\mathcal{E}$ step by constructing a fourth neural network
$g : \Sigma \times [0,T] \times \Theta_g \rightarrow \widehat{\Gamma}( \text{Sym}^2 T^* \mathcal{M}),   \widehat{\Gamma} \subseteq \mathcal{F}_{\epsilon}(\Gamma) = \Big\{ \widetilde{g} \in \Gamma(\text{Sym}^2 T^* \mathcal{M}) \ |  \ || \widetilde{g} - g_{\theta_g} || < \epsilon \Big\},$
where the second equation is used to denote that our simulated Riemannian is near the true Riemannian metric in the function space in an $\epsilon$-neighborhood, denoted $\mathcal{F}_{\epsilon}$. We have used notation $\widehat{\Gamma}( \text{Sym}^2 T^* \mathcal{M})$ to illustrate the space of approximate smooth symmetric 2-tensors on the tangent space of the manifold. In particular, our neural network approximates the true Riemannian metric such that $g_{\theta_g,ij}(p) \approx g_{ij}(p) = g_p \Big( \frac{\partial}{\partial u^i} \Big|_p, \frac{\partial}{\partial u^j} \Big|_p \Big), g \in \Gamma( \text{Sym}^2 T^* \mathcal{M}) $\
for tangent basis vector $\partial / \partial u^i$.

\vspace{2mm}

\begin{wrapfigure}{r}{0.35\textwidth}
  \centering
  \includegraphics[width=\linewidth]{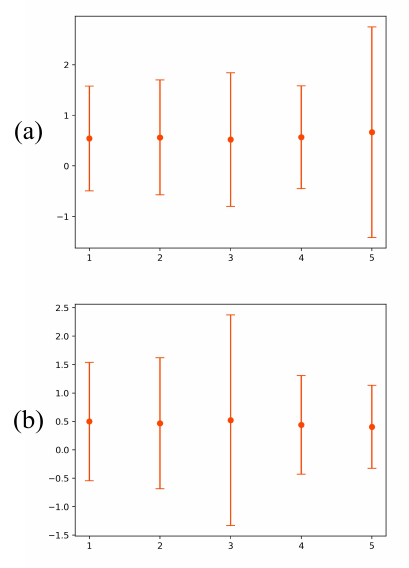}
  \caption{We plot the means and standard deviations of errors on a specific out-of-distribution setting with a (a) vanilla VAE and (b) the Perelman's functional experiment. With our hyperparameters, results are consistent, especially means.}
  \label{fig:errors_consistent}
  \vspace{-40mm}
\end{wrapfigure}

We propose a physics-informed neural network (PINN) framework \cite{wang2020understandingmitigatinggradientpathologies} \cite{wang2023expertsguidetrainingphysicsinformed} \cite{raissi2017physicsinformeddeeplearning} for finding the metric solution $g$ to the geometric flow. In a tradtional objective, the $L^2$ integral loss can be approximated via Monte-Carlo integration in a discrete-type objective $\EX_{t \sim U([0,T])} \Big[ \EX_{u \sim \rho_u/\int \rho_u du} \Big[ \Big|\Big| \partial_t g_{\theta_g}(u,t) + \mathcal{F}_K[g_{\theta_g}(u,t)] \Big|\Big|_F^2  \Big] \Big] $. Our approaches will be consistent with this general idea but partially deviate from this. Typically, in a PINN, the physics loss is evaluated uniformly with respect to its domain: this is not the case for us, and this is where $\rho$ serves its purpose. $\rho$ is a weighting function so that the physics loss is uniform with respect to input $\tilde{\phi}_0$ into $u$. In particular, we have $
\rho_u = (u_{\theta_u})_\# \mu_{\Phi} = (u_{\theta_u})_\# \mathcal{U}_{\Phi},  (u_{\theta_u})_{\#} :  \mathcal{U}_{\Phi_{\text{supp}(\mathcal{X})}} \rightarrow \mathcal{M}^+(\Sigma)$
as a pushforward of a measure between spaces equipped with the Borel sigma algebras $\sigma$.  Our primary losses will be in a typical encoder-decoder fashion with loss of the form
\begin{gather}
\inf_{\Theta} \int \Bigg[ \EX_{ u \sim \tilde{q}} \Big[ || \mathcal{D}_{\theta_{\mathcal{D}}} \circ \mathcal{E}_{\theta_{\mathcal{D}}} \circ u_{\theta_u} ( \tilde{\phi}_0, t) - \tilde{\phi}_t ||_2^2 \Big] 
\\
+ \frac{\beta}{T} \EX_{\tilde{\phi}_0} \Big[ D_{KL} ( \tilde{q}(u |  \tilde{\phi}_0) \| p(u) ) \Big] 
\\
+ \EX_{\tilde{\phi}_0} \Big[ || \text{geometric flow PDE} ||_F^2 + ||g - J^T J||_F^2 \Big] \Bigg] dt  .
\end{gather}
Here $J$ is the Jacobian of the parametric map $\mathcal{E}$ with respect to local coordinates. Hence, $J^T J$ is the induced metric. We will  consider variational settings \cite{kingma2022autoencodingvariationalbayes} \cite{Kingma_2019} \cite{Girin_2021} with KL divergence regularizing loss in latent spaces.

\section{Contributions}

\subsection{Ricci flow on the sphere}

Curvature is already uniformized as the sphere, and Ricci flow maintains the fact that the geometric structure over the entire flow is the same, except under Riemannian metric rescalings. In this case, we will normalize our data and rescale it under the correct radius subject to Ricci flow up to constants using $
\mathcal{E}_{\mathbb{S}^d} = R(t) \times \frac{u}{||u||_2} = \sqrt{ R^2 - C(d-1)t} \times \frac{u}{||u||_2} $.

\subsection{Path integration with Gaussian curvature} 

We derive a special case for 2-manifolds using variations of Gaussian curvature \cite{elemofdiffgeo} \cite{weisstein_gaussian} under restrictions, as this allows closed path integration upon the curvature. Note that Ricci flow for 2-manifolds can be expressed using Gaussian curvature $K$, thus this section is exactly Ricci flow for 2-manifolds. Our geometric flow loss is motivated by the classical PDE $
(\partial_t - K) g = 0 $.
Our total geometric flow loss function in this setting is
\begin{align}
\int \Bigg( \EX_{\rho_{\Sigma}} \Bigg[ \sqrt{\text{det}(g)} \sum_{ij} \frac{\Big(  \partial_t g \circ \partial_t g  - \widetilde{K}^2  g \circ g  \Big)_{ij}}{  \Big( \partial_t g \circ g - 2 \widetilde{K} g \circ g \Big)_{ij}} \Bigg]  +  \EX_{(u_0,r) \sim \rho_{\Sigma} \otimes \rho_r} \Bigg[\oint \frac{\sqrt{\text{det}(g)}}{g_{11}}({\Gamma_{11}}^2 du^1  + {\Gamma_{12}}^2 du^2 ) \Bigg] \Bigg) dt .
\end{align}
The contour integral makes computation  efficient; we can sample any ball, and the highest order computation needed is only the first order Christoffel symbols. $\circ$ is the Hadamard (element-wise) product. $\widetilde{K}$ is an approximate surrogate for a Gaussian curvature estimation.

\subsection{A new flow derived with a Taylor expansion and low-order curvature proxies}
\label{sec:fundamental_form}

Our motivation for this flow is that Ricci curvature is expensive to compute, and the Weyl tensor does not have a closed form, so we replace these with properly justified proxies. We develop a custom geometric flow based on the Riemannian tensor decomposition with the Weyl tensor omitted. In particular, the flow of this next section is not as much derived as it is invented, motivated by the Riemannian decomposition and a linear Gauss equation. We truncate Riemannian curvature and replace the Weyl tensor with a new tensor (see Appendix \ref{app:second_FF_sec}) $
R_{ijkl} = W_{ijkl} + S_{ijkl} + U_{ijkl} \rightarrow\widetilde{R}_{ijkl} =  S_{ijkl} + U_{ijkl}$ since the Weyl tensor is trace-free.
We examine the Gauss equation 
\cite{Gauss_equation} 
\begin{align}
\langle R(X,Y)Z, W \rangle = \langle \overline{R}(X,Y)Z,W \rangle + \langle \Pi (X, W), \Pi(Y,Z) \rangle  - \langle \Pi(X,Z), \Pi(Y,W) \rangle ,
\end{align}
in the setting that $\overline{R}_{ijkl} = 0$. Notice that this equation is quadratic in the second fundamental form $h$ by the Gauss equation, but we will draw a connection from this equation but linearized. We will operate so that a Taylor approximation at $\Pi_0$ is nearly satisfied,
\begin{align}
f(\Pi) = f(\Pi_0) + \Bigg( \Pi_{ik,0} \delta \Pi_{jl} + \delta \Pi_{ik} \Pi_{jl,0} - \Pi_{jk,0} \delta \Pi_{il} - \delta \Pi_{jk} \Pi_{il,0} \Bigg) + \text{higher order terms} .
\end{align}
This linearization is especially motivated in the case of the sphere (see Appendix \ref{app:second_FF_sec}), which we use in our experiments (see section \ref{sec:burgers_experiments}). This expansion is meaningful because under a uniform-curvature metric and in our framework
\begin{align}
\text{Ric}_{ik} = \kappa \Bigg[ (d-1) \Pi_{ik,0}  +  g_{ik,0} \text{Tr}(\delta \Pi) + (d-2)\delta \Pi_{ik} \Bigg] + \mathcal{O}(\delta \Pi^2) ,
\end{align}
thus $\text{Ric}_{ij}$ and $\Pi_{ij}$ scale similarly, meaning we can replace Ricci curvature with a low-order surrogate. We will set the traceless extrinsic curvature tensor $
\widehat{H}_{ij} = \Pi_{ij} - \frac{H}{d} g_{ij} $.
We will compute $h$ with a Hessian surrogate in the normal bundle $
\Pi_{ij} = | | \partial_{ij} \mathcal{E} - \sum_{kl} g^{kl} \langle \partial_{ij} \mathcal{E}, \partial_k \mathcal{E} \rangle \partial_l \mathcal{E} | |_2 $. Our geometric flow loss is $
\EX \Big[ \ || \partial_t g_{ij} - \text{constant} \cdot H_{ij} ||_F^2  \ \Big] $, where $H_{ij}$ is defined as in Appendix \ref{app:second_FF_sec}. This flow is motivated by the fact that, under a spherical geometry, Taylor-expanded linear Gauss equation approximates the true quadratic equation.

\subsection{Perelman's functional and non-parametric geometric flows with scalar curvature special cases}

Our approaches are similar to \cite{lazarev2023metric} here. We propose considering the conformally changed metric $g(u,t) =  \psi^{4/(d-2)}(u,t) \cdot g $
with special case of scalar curvature
$R(\psi^{4/(d-2)}g) = \frac{ 4 \frac{d-1}{d-2} \Delta_g \psi - R(g) \psi}{ \psi^{ (d+2)/(d-2)}} $. We derive a result based on Perelman's $\mathcal{W}$ \cite{perelman2002entropyformularicciflow} functional and conformally changed metrics. We also tie connections between Monte-Carlo integration and physics-informed methods, since their loss is of equivalence. By inducing a nonzero time derivative on this functional, a geometric flow is induced. This result is nonparametric because the geometric flow is not easily quantifiable, but is rather held using the data with neural networks. The loss we develop here is
\begin{align}
\Bigg| \text{positive constant} -  \frac{d}{dt} \int_{\Sigma}  \frac{ 4 \frac{d-1}{d-2} \Delta_g \psi - R(g) \psi}{ \psi^{ (d+2)/(d-2)}} \sqrt{\text{det}(g(u,t))} \rho(u) \bigwedge_k du^k \Bigg| .
\end{align}

\subsection{Geometric flows with differentiating harmonic map energy}

In this section, we work with extrinsic harmonic maps \cite{HeleinHarmonic} \cite{Chervon_2004}. We construct a harmonic map with a neural network $
\psi_{\theta_{\psi}} : \mathcal{M}_0 \hookrightarrow \mathcal{N}_t \subset \mathbb{R}^D $
that is defined \textit{extrinsically} so that it is embedded in the space. We develop a loss based on the energy $
E(\psi)(t) = \int_{\Sigma} |d \psi|^2 du$
and its derivative so that it is sufficiently nonzero. Using the product rule under regularity conditions (see Appendix \ref{app:harmonic_maps}), we are led to the loss function
\begin{align}
\label{eqn:harmonic_energy_loss}
\Bigg| \text{constant} - & \Bigg|  \int_{\Sigma} g^{ij}(u) \Big( \Big\langle  \frac{\partial}{\partial t} \frac{\partial \psi(\varphi(u))}{\partial u^i}, \frac{\partial \psi(\varphi(u))}{\partial u^j} \Big\rangle_{\mathbb{R}^D}   + \Big\langle  \frac{\partial \psi(\varphi(u))}{\partial u^i}, \frac{\partial}{\partial t} \frac{\partial \psi(\varphi(u))}{\partial u^j} \Big\rangle_{\mathbb{R}^D}  \Big) \sqrt{\text{det}(g)} \bigwedge_k du^k  \Bigg|  \Bigg| .
\end{align}

\vspace{2mm}

\SetKwComment{Comment}{/* }{ */}

\begin{table}[!htbp]
\caption{\textcolor{Cerulean}{Variational} method errors by scenario with the means and standard deviations listed using an $L^1$ relative error metric. We illustrate our methods over 300 out-of-distributions data samples.  Our OOD scenarios are: (1) noise with $\sigma=1.0$ injected in all 201 locations; (2) noise with $\sigma=1.5$ injected in all 201 locations; (3) noise with $\sigma=2.5$ injected in all 201 locations; (4) noise with $\sigma=2.0$ injected in 31 locations; (5) noise with $\sigma = 3.0$ injected in 11 locations; (6) noise with $\sigma=0.25$ injected in all 201 locations; (7) initial data scaled by 1.4 with noise of $\sigma=1.5$ injected in 101 locations; (8) a linear function of the form $y_0 + 2 \pi j\Delta x/N$ injected everywhere; (9) noise with $\sigma=3.0$ injected in the first 11 locations only. As a last remark, we found our VAE baseline was generally pretty consistent among retraining (see Figure \ref{fig:errors_consistent}). For example, we retrained the baseline VAE and repeated scenario 1 with a $0.727$ mean.}
\label{tab:method_errors}
\centering
\renewcommand{\arraystretch}{1.2}%
\footnotesize
\begin{tabular}{
  >{\raggedright\arraybackslash}p{2.5cm}
  >{\centering\arraybackslash}p{2.25cm}
  >{\centering\arraybackslash}p{2.25cm}
  >{\centering\arraybackslash}p{2.25cm}
  >{\centering\arraybackslash}p{2.25cm}
  >{\centering\arraybackslash}p{2.25cm}
}
\toprule
\multirow{2}{*}{Scenario} & \multicolumn{5}{c}{Method [$\downarrow$]} \\
\cmidrule(lr){2-6}
& Path integration & Second fundamental form flow & Perelman's functional & Harmonic functional & Baseline \\
\midrule
\rowcolor{Emerald!15} Scenario 1 & {$0.592 \pm 1.270$} & {$0.555 \pm 1.512$} & {$\textBF{0.482 \pm 0.962}$} & {$0.522 \pm 1.244$} & {$0.765 \pm 2.704$} \\
\addlinespace
\rowcolor{LimeGreen!15} Scenario 2 & {$0.736 \pm 1.498$} & {$0.735 \pm 1.986$} & {$\textBF{0.643 \pm 1.188}$}  & {$0.699 \pm 1.656$} & {$0.984 \pm 3.217$} \\
\addlinespace
\rowcolor{Emerald!15} Scenario 3 & {$0.944 \pm 1.714$} & {$0.961 \pm 2.587$} & {$\textBF{0.853 \pm 1.376}$} & {$1.032 \pm 2.317$} & {$1.231 \pm 3.442$} \\
\addlinespace
\rowcolor{LimeGreen!15} Scenario 4 & {$0.377 \pm 0.771$} & {$\textBF{0.361 \pm 0.720}$} & {$0.392 \pm 1.504$} & {$0.389 \pm 1.248$} & {$0.487 \pm 1.971$} \\
\addlinespace
\rowcolor{Emerald!15} Scenario 5 & {$0.339 \pm 0.902$} & {$0.343 \pm 0.724$} & {$0.355 \pm 1.432$} & {$\textBF{0.334 \pm 0.799}$} & {$0.408 \pm 1.482$} \\
\addlinespace
\rowcolor{LimeGreen!15} Scenario 6 & {$0.273 \pm 0.514$} & {$0.182 \pm 0.510$} & {$\textBF{0.152 \pm 0.302}$} & {$0.245 \pm 0.356$} & {$0.237 \pm 0.881$} \\
\addlinespace
\rowcolor{Emerald!15} Scenario 7 & {$0.354 \pm 0.621$} & {$\textBF{0.348 \pm 0.857}$} & {$0.408 \pm 1.146$} & {$0.387 \pm 0.652$} & {$0.523 \pm 2.263$} \\
\addlinespace
\rowcolor{LimeGreen!15} Scenario 8 & {$\textBF{0.610 \pm 0.377}$} & {$1.113 \pm 1.948$}  & {$0.613 \pm 0.331$} & {$0.727 \pm 0.257$} & {$1.310 \pm 2.675$} \\
\addlinespace
\rowcolor{Emerald!15} Scenario 9 & {$0.432 \pm 1.120$} & {$0.420 \pm 0.839$} & {$\textBF{0.416 \pm 0.912}$} & {$0.417 \pm 0.881$} & {$0.468 \pm 1.491$} \\
\addlinespace
\bottomrule
\end{tabular}
\end{table}

\bibliographystyle{plainnat}
\bibliography{bibliography}

\appendix

\newpage

\section{Additional discussion}

Our overarching arguments are established so that we attempt to present empirical justification that different geometric representations have different ambient outcomes in a beneficial way: structured representations facilitate robustness. This argument is counterposed to that in which we argue manifold representations themselves aid robustness. Indeed, the phenomenon that latent spaces can, but not necessarily always, naturally acquire Riemannian geometric structure, notably manifolds, on their own accord has been observed in certain settings, primarily \cite{chadebec2022geometricperspectivevariationalautoencoders}. This work specifically differs from ours in the senses of low-dimensionality and controllability. \cite{chadebec2022geometricperspectivevariationalautoencoders} argues more closely that taking into account the natural geometry that exists can improve performance in areas such as sampling. Our work is more closely aligned with regulating these geometries and their low-dimensional structures.

\vspace{2mm}

We argue the parameterized geometric flow is a regularizer. One approach to regularize a latent space in a manifold-based approach is through projections, in which \cite{lopez2025gdvaesgeometricdynamicvariational}\cite{lopez2021variationalautoencoderslearningnonlinear} study this aforementioned technique. Our approaches primarily investigate when latent spaces are imposed with a manifold structure upon in-distribution data such that the data naturally formulates this manifold. Thus, we maintain faithful geometric representations across the entire time frame of the ambient data, perfect for PDE-type data. This notion is consistent with the mesh-invariance archetype of neural operators \cite{Lu_2021}, and our methods are aligned with mesh invariance over the time domain.

\vspace{2mm}

Curvature-type flows \cite{toppingricciflow} \cite{perelman2002entropyformularicciflow} \cite{ricciflowsheridan} \cite{calegariricciflow} will be our means of operation, but we will develop new intrinsic geometric flows derivative of curvature flows that are relevant to the beneficial computational properties that Ricci flow brings forth. The following are the theoretical properties of curvature flows that will aid us in our empirical studies:
\begin{itemize}
    \item It is not true that $||g||_F \rightarrow 0$ implies that $\text{curvature} \rightarrow 0$. Since our flows are curvature-based, this will imply $\partial_t g$ is sufficiently large over the entire course of the flow. This property helps us in that a near-triviality metric is not learned in all settings. For example, if we take $g$ to be small, then the entire flow and the overall manifold's measure would also remain small under a compact time interval of Ricci flow. Since this is not true, this helps ensure nondegeneracy. Consider the example
    \begin{align}
    \tilde{g} = \epsilon^2 g_{\mathbb{S}^d}, \text{Ric}(\tilde{g}) = \text{Ric}(\epsilon^2 g_{\mathbb{S}^d} ) = \text{Ric}(g_{\mathbb{S}^d}) = (d-1)g_{\mathbb{S}^d} =  \frac{(d-1)}{\epsilon^2} \tilde{g} .
    \end{align}
    We have denoted $\mathbb{S}^d$ the $d$-sphere. Therefore, Ricci curvature is bounded below as the metric extinguishes (in fact it diverges). We have observed that an arbitrarily small metric is antithetical towards robustness, thus this property is valuable.
    \item The volume of the manifold expands or contracts heavily under curvature flows. Significant manifold deformation helps the learning task. In particular, the volume element under Ricci flow evolves according to
    \begin{equation}
    \frac{d}{dt} \sqrt{\text{det}(g)} = - \text{Ric}(g)_{ij} g^{ij} \sqrt{\text{det}(g)} .
    \end{equation}
    This result comes from Jacobi's formula. We also remark $\text{Ric}_{ij} g^{ij}$ is precisely scalar curvature. We require $\text{Ric}_{ij} g^{ij} \neq 0$. This property is valuable because a distortive geometry helps the learning task. A constant geometry means two distinct ambient times cannot be discerned as distinct in the latent space. We note that, with intrinsic flows, our geometries can displace in an immersion. By volume element distortion, there is greater representation capacity and flexibility. Latent representations can be conveyed by both immersion displacement and deformation itself. In our experiments, we found a geometry only displacing independent of the flow does not yield sufficiently good inference quality. We discovered this through our second fundamental form experiments without the use of our $\Gamma_{ijkl}$ diagonal term (see Appendix \ref{app:second_FF_sec}).
\end{itemize}
We make some other remarks. Extrinsic flows would require differential equations on vectors orthogonal to the tangent space, up to scaling such as via mean curvature. To simulate such an evolution, one strategy would be to employ the use of neural ordinary differential equations. Aside from computing necessary normal vectors and scaling coefficients, neural ODEs are nontrivial in large scale settings. With our methods, computing tangent vectors is cheap, and simulating a Riemannian metric, aside from Ricci curvature computation, is generally easier. Hence, we have empirically found intrinsic flows work nicely. Our intrinsic passively geometric flows exert their own vector fields in our latent spaces. Curvature-based flows create their own manifold vector field, such as through vector fields
\begin{gather}
\xi = \mathcal{K}^{\sharp}(X)  \in \Gamma(T \mathcal{M}), \ \ \ \ \ \mathcal{K}^{\sharp} \in \Gamma(\text{End}(T \mathcal{M}))
\\
\xi' = \Big( \mathcal{K}^{\sharp}(X) \Big)^{\perp}  \in \Gamma(N \mathcal{M})
\end{gather}
for correct $X$ in our latent immersions. The musical isomorphism $\sharp : T^* \mathcal{M} \xrightarrow{\cong} T \mathcal{M}$ maps the respective tensor from the cotangent bundle to a vector field in the tangent bundle. $\mathcal{K}^{\sharp}$ is a generalized vector field built on some idea of curvature, i.e. $\text{Ric}^{\sharp}$. $\Gamma$ is the space of smooth vector fields on the manifold, and $\text{End}$ denotes an endomorphism with the tangent bundle. Thus, our intrinsic flows behave similarly to extrinsic in certain ways.

\vspace{2mm}

Our methods pair with the manifold hypothesis, but this rather states datasets themselves lie on intrinsically-low-structured manifolds in a high dimensional space, not that the latent spaces themselves formulate into manifolds. Again, we reference \cite{chadebec2022geometricperspectivevariationalautoencoders}, and our arguments are not tailored such that whether or not geometries exist in these latent spaces, but rather that we can regulate the structure that developments to induce outcome.

\section{Baseline training}

\textbf{\begin{figure}[!h]
  \vspace{0mm}
  \centering
  \includegraphics[scale=0.45]{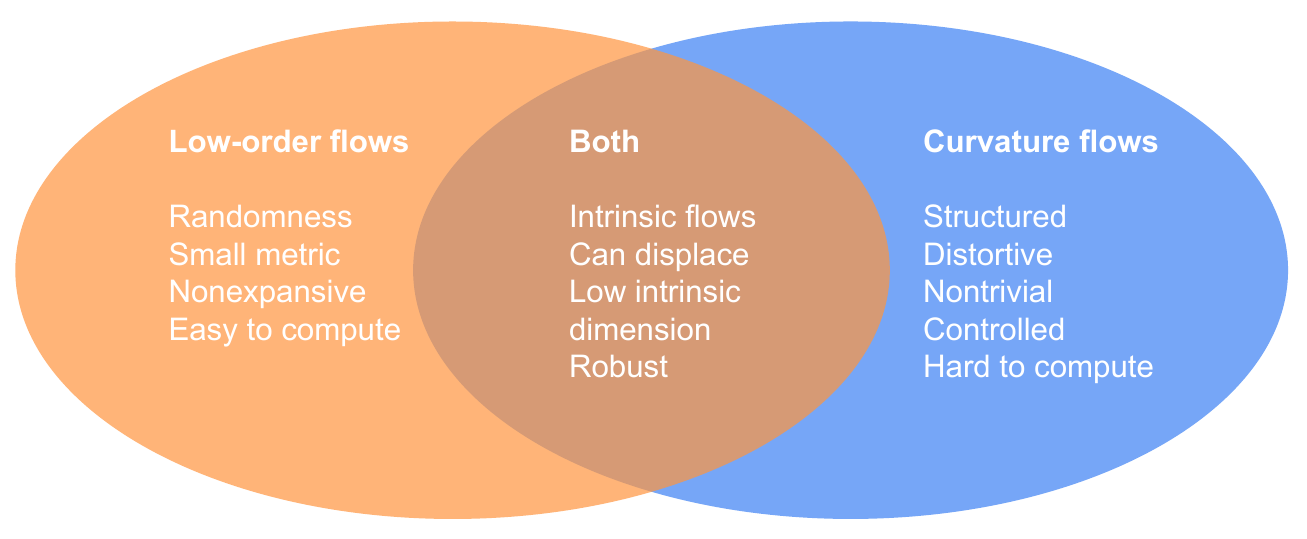}
  \caption{We compare curvature-type flows to alternative flows that are easy to compute, which are typically trivial in derivative order.}
\label{fig:flow_venn_diagram}
\end{figure}}

\begin{algorithm}[hbt!]
\caption{Geometric flow-regularized encoder-decoder training}\label{alg:two}
\KwData{$A_1 = \{ \tilde{\phi}_0^{(i)} \}_i, A_2 = \{ \tilde{\phi}_{t_j}^{(i)} \}_{ij}$ evaluated at $(x_1,\hdots,x_n,t_j) \in \mathcal{X} \times [0,T]$, $u_{\theta_u}, \mathcal{E}_{\theta_{\mathcal{E}}}, \mathcal{D}_{\theta_{\mathcal{D}}}$ }
\While{\text{Loss not converged}}{
 Sample $N$, $A_1, A_2, t$\;
 Evaluate $u(\tilde{\phi}_0^{(i)})$\;
 Map to manifold $\mathcal{E}(u)$\;
 Evaluate $|| g(u,t) + \mathcal{F}_K[g(u,t)]||_F^2$\;
 Compute $\big\langle \partial_k \mathcal{E}(u,t), \partial_l \mathcal{E}(u,t) \big\rangle$\;
 Decode point on manifold $\tilde{d} = \mathcal{D} \circ \mathcal{E}$\;
 Evaluate $||\tilde{d} - \tilde{\phi}_{t} ||_2^2$\;
 Formulate total objective and take grad\;
}
\end{algorithm}

\begin{table}[!htbp]
\caption{Training time by method (in seconds). We have used intrinsic dimension $d=3$ for all experiments except the Ricci flow ($d=2$), and for the Gaussian curvature (since Gaussian curvature only exists for $d=2$). We use 100 descent iterations.}
\label{tab:method_errors}
\centering
\renewcommand{\arraystretch}{1.2}%
\begin{tabular}{
  >{\raggedright\arraybackslash}p{1.8cm}
  >{\centering\arraybackslash}p{1.8cm}
  >{\centering\arraybackslash}p{1.8cm}
  >{\centering\arraybackslash}p{1.8cm}
  >{\centering\arraybackslash}p{1.8cm}
  >{\centering\arraybackslash}p{1.8cm}
  >{\centering\arraybackslash}p{1.8cm}
}
\toprule
\multirow{2}{*}{} & \multicolumn{6}{c}{Method} \\
\cmidrule(lr){2-7}
& Path integration & Second fundamental form flow & Perelman's functional & Harmonic functional & Ricci flow ($d=2)$ & Baseline (extended VAE) \\
\midrule
\rowcolor{Dandelion!15} Training time & 10.115 & 18.602 & 11.657 & 10.348  & 27.005 & 6.162 \\
\addlinespace
\bottomrule
\end{tabular}
\end{table}

\section{Setup}

We will be learning geometric flows of a generalized form
\begin{equation}
\partial_t g + \mathcal{F}_{K}[g] = 0  .   
\end{equation}
Here, we have denoted $\mathcal{F}$ a generalized differential operator that is relevant to curvature, denoted $K$, i.e.
\begin{align}
\mathcal{F}_K[g] = \mathcal{A}_K g, \ \ \ \ \ \mathcal{A}_K \in \mathcal{A}_{\text{nonlin},\text{lin}}(L^2(\mathcal{M},g))  . 
\end{align}
We have used $\mathcal{A}$ as notation for the set of nonlinear and linear bounded operators on Riemannian metric $g$,
\begin{align}
& \mathcal{A}_{\text{nonlin},\text{lin}}(L^2(\mathcal{M},g)) =  \Big\{ \mathcal{A} : L^2(\mathcal{M}, \text{Sym}^2 (T^* \mathcal{M}) ) \rightarrow  L^2(\mathcal{M}, \text{Sym}^2 (T^* \mathcal{M}) ) : \mathcal{A} \ \text{is bounded and linear} \Big\} 
\\
& \ \ \ \ \ \ \ \ \ \ \ \ \ \ \ \ \ \cup \ \Big\{ \mathcal{A} : L^2(\mathcal{M}, \text{Sym}^2 (T^* \mathcal{M}) )\rightarrow   L^2(\mathcal{M}, \text{Sym}^2 (T^* \mathcal{M}) ) : \mathcal{A} \ \text{is bounded and nonlinear} \Big\}   .
\end{align}
We have not assumed our manifold is compact in this theoretical setting (in fact, for us, it will be a discrete point cloud immersed in a space with respect to training data). There are many variations of curvature, so $K$ is just a generalized notion of this. We will specifically be studying variations of Ricci curvature, but not Ricci curvature exactly. 

\vspace{2mm}

For our ambient data, we will restrict our attention to PDEs whose image is a subset of $\mathbb{R}$, i.e. $\phi_t \in \mathbb{R}$. We seek a mapping
\begin{gather}
\Gamma : (\phi_0, t) \in C(\mathcal{X} ;\mathbb{R}) \cap W^{k,p}(\mathcal{X}; \mathbb{R}) \times [0,T] \rightarrow \phi_t \in C(\mathcal{X} ; \mathbb{R}) \cap W^{k,p}(\mathcal{X}; \mathbb{R}), 
\\
\text{or discretized,} 
\\
\Gamma^{\dagger} :  \Phi \subseteq \Bigg\{ \tilde{\phi}_0 \in \mathbb{R}^N : \tilde{\phi}_0 = \phi_0 |_{\Omega}, \phi_0 \in C(\mathcal{X}; \mathbb{R}) \cap W^{k,p}(\mathcal{X}; \mathbb{R}) \Bigg\} \times [0,T]  \\
 \ \ \ \ \ \ \ \ \ \ \ \ \ \ \ \ \ \ \ \ \  \longrightarrow \Phi' \subseteq \Bigg\{ \tilde{\phi}_t \in \mathbb{R}^N : \tilde{\phi}_t = \phi_t |_{\Omega}, \phi_t \in C(\mathcal{X}; \mathbb{R}) \cap W^{k,p}(\mathcal{X}; \mathbb{R}) \Bigg\}.
\end{gather}
Here, $k$ and $p$ are suitable choices of constants, and we have discretized the dynamics as vectors in $\mathbb{R}^N$ and evaluated along the mesh in space and time
\begin{align}
\Omega_{\Delta x} \times \mathcal{T}_{\Delta t} = \Big\{ (x_0 + k \Delta x, j \Delta t) : k \in \mathbb{N}^{\text{dim}(\mathcal{X})}, j \in \mathbb{N}, (\max_j j) \Delta t = T \Big\} \subset \mathcal{X} \times \mathbb{R}^+ .
\end{align}
We will take this mesh to be equispaced.

\vspace{2mm}

We will employ the Whitney immersion theorem in all of our experiments, and so all of our $d-$manifolds will be immersed in $2d-1$ Euclidean space unless otherwise stated, which are primarily are our Gaussian curvature experiment, which is restricted to surfaces, and the sphere experiment.

\vspace{2mm}

We will assume every manifold is orientable, thus we can integrate against the wedge product.

\section{A vanishing Ricci off-diagonal: geometric meaning versus computational gains}

In our experiments, it would notable to consider restricting both the metric and the Ricci tensor as diagonal for computational expense purposes. We provide the following result that this is not always meaningful. In particular, we show just because the metric is diagonal at some $t \in [0,T]$, it is not necessarily diagonal for all $\tau>t$ due to the Ricci tensor. 

\vspace{2mm}

\textbf{Theorem 1.} Let $g, \text{Ric} \in \Gamma(T^* \mathcal{M} \otimes T^* \mathcal{M})$ be sufficiently smooth on $[0,T], T > 0$, and let $g$ be diagonal at $t=0$. Let $g$ be subject to Ricci flow on $[0,T]$. Then $\text{Ric}$ is not necessarily diagonal for some $t$, hence neither is $g$, over $[0,T]$.

\vspace{2mm}

\textit{Proof.} Let $i \neq k$. We will use Einstein notation only if an index is repeated in lower and upper positions exactly once each. We begin by noting Hamilton's result of the time evolution of the Ricci tensor under Ricci flow \cite{hamilton1982three}
\begin{equation}
\label{eqn:ricci_tensor_time_deriv}
\partial_t \text{Ric}_{ik} = \Delta \text{Ric}_{ik} + 2g^{pr} g^{qs} R_{piqk} \text{Ric}_{rs} - 2g^{pq} \text{Ric}_{pi} \text{Ric}_{qk} .
\end{equation} 
Suppose necessarily that $g_{ij} = \text{Ric}_{ij} = 0$ as well as their spatial derivatives at $t=0$, i.e. the Ricci tensor is identically $0$ along the off-diagonal. We will show $\Delta \text{Ric}_{ik}$ is not necessarily $0$ due to contractions over diagonal terms, nor is $\partial_t \text{Ric}_{ik}$. First, we compute $\Delta \text{Ric}_{ik}$. First, note if the metric is diagonal, so is its inverse. Observe using the covariant derivative formulas (see \ref{eqn:cov_deriv_formula_2}, \ref{eqn:cov_deriv_formula_3})
\begin{align}
& \Delta \text{Ric}_{ik} = g^{jl} \nabla_j \nabla_l \text{Ric}_{ik} = g^{jl} \nabla_j ( \partial_l \text{Ric}_{ik} - {\Gamma_{il}}^m \text{Ric}_{mk} - {\Gamma_{kl}}^m \text{Ric}_{im} )  
\\[1em]
& = \sum_j g^{jj} \Big[ \partial_j (\nabla_j \text{Ric}_{ik}) - {\Gamma_{jj}}^m \nabla_m \text{Ric}_{ik} - {\Gamma_{ji}}^m \nabla_j \text{Ric}_{mk} - {\Gamma_{jk}}^m \nabla_j \text{Ric}_{im} \Big]
\\[1em]
&= \sum_j \Bigg[ g^{jj} \Bigg( \partial_j  ( \partial_j \text{Ric}_{ik} - {\Gamma_{ij}}^m \text{Ric}_{mk}  - {\Gamma_{kj}}^m \text{Ric}_{im} )  - {\Gamma_{jj}}^m ( \partial_m \text{Ric}_{ik} - {\Gamma_{im}}^k \text{Ric}_{kk} - {\Gamma_{km}}^i \text{Ric}_{ii} ) 
\\
& \ \ \ \ \ \ \ \ \ \ \ \ - \sum_m  {\Gamma_{ji}}^m (  \partial_j \text{Ric}_{mk} - {\Gamma_{mj}}^k \text{Ric}_{kk} - {\Gamma_{kj}}^m \text{Ric}_{mm})  - \sum_m {\Gamma_{jk}}^m ( \partial_j \text{Ric}_{im} - {\Gamma_{ij}}^m \text{Ric}_{mm} -  {\Gamma_{mj}}^i \text{Ric}_{ii})  \Bigg) \Bigg]  .
\end{align}
Let us now apply our diagonal constraints. Note the symmetry in the expression $- {\Gamma_{ij}}^k \text{Ric}_{kk} - {\Gamma_{kj}}^i \text{Ric}_{ii}$. Now, we note the following: ${\Gamma_{ij}}^k = 0$ for diagonal metrics if $i \neq j \neq k$ \cite{diagonal_christoffel_symbols}. For the first term, $j=i,k$ otherwise we get triviality. For the inner summations, we must have $m = j,i,k$, and $j$ can run over any index. Hence,
\begin{align}
= & \sum_{j \in \{i,k\}} \Bigg[ g^{jj} \partial_j ( - {\Gamma_{ij}}^k \text{Ric}_{kk} - {\Gamma_{kj}}^i \text{Ric}_{ii} ) \Bigg]
\\ & + \sum_j g^{jj} \sum_m {\Gamma_{jj}}^m \Bigg[ {\Gamma_{im}}^k \text{Ric}_{kk} + {\Gamma_{km}}^i \text{Ric}_{ii} \Bigg]
\\ & + \sum_j g^{jj} \Bigg[ - {\Gamma_{ji}}^k \partial_j \text{Ric}_{kk} + \sum_{ m \in \{j,i,k\} \ \text{distinct}} \Big( {\Gamma_{ji}}^m {\Gamma_{mj}}^k \text{Ric}_{kk} + {\Gamma_{ji}}^m{\Gamma_{kj}}^m \text{Ric}_{mm} \Big) \Bigg]
\\ & + \sum_j g^{jj} \Bigg[ - {\Gamma_{jk}}^i \partial_j \text{Ric}_{ii} + \sum_{ m \in \{j,i,k\} \ \text{distinct}} \Big( {\Gamma_{jk}}^m {\Gamma_{ij}}^m \text{Ric}_{mm} + {\Gamma_{jk}}^m {\Gamma_{mj}}^i \text{Ric}_{ii} \Big) \Bigg] := (a) .
\end{align}
Now, returning to equation \ref{eqn:ricci_tensor_time_deriv}, we have
\begin{align}
\label{eqn:eqn_38}
\partial_t \text{Ric}_{ik} & = (a) + \sum_{p,q} 2 g^{pp} g^{qq} R_{piqk} \text{Ric}_{pq} - 2 g^{ik} \text{Ric}_{ii} \text{Ric}_{kk} = (a) + \sum_q 2 (g^{qq})^2 g_{qm} {{R}_k}{}^{m}{}_{qi} \text{Ric}_{qq} - 2 g^{ik} \text{Ric}_{ii} \text{Ric}_{kk}  \displaybreak[1]
\\ & = (a) + \sum_q 2 \Bigg[  (g^{qq})^2 g_{qq}  {{R}_k}{}^{q}{}_{qi}
\text{Ric}_{qq}  \Bigg] - 0 = (a) + \sum_q 2 \Bigg[ \underbrace{ g^{qq} {{R}_k}{}^{q}{}_{qi}
\text{Ric}_{qq} }_{=(b)}  \Bigg] - 0 = (c) .
\end{align}
We remark the diagonalization of $g$ reduces the summation over $m$ in $\ref{eqn:eqn_38}$ to a single term, which means that ${{R}_k}{}^{q}{}_{qi}$ is not Einstein notation and does not reduce to a Ricci curvature term. In general, $(c)$ is not identically $0$. Observe it is possible terms involving a product of the Riemannian tensor and Ricci curvature appear in $(a)$ and $(b)$, but we need not do this computation. The first summation in $(a)$ only runs over $j=i,k$, but the summation in $(b)$ runs over any index. This is a guarantee that not all terms necessarily cancel out (specifically for intrinsic dimension $> 2$). We provide a non-rigorous mathematical argument why the result follows in the remark below.

\vspace{2mm}

Thus, the time derivative is nonzero, so the off-diagonal of $\text{Ric}$ is not identically zero and we are done. In particular, we have
\allowdisplaybreaks
\begin{align}
g_{ij}(t) & = g_{ij}(0) -2 \int_0^t (\text{Ric}_{ij}(0) + \int_0^{\tau} \partial_{t} \text{Ric}_{ij}(\gamma) d\gamma ) d\tau    = -2 \int_0^t \int_0^{\tau} \partial_{t} \text{Ric}_{ij}(\gamma) d\gamma d\tau     
\end{align}
is not necessarily 0. We remark we also use smoothness, and so $\partial_t \text{Ric}_{ik}$ is not only nonzero on a set of measure 0.
$\square$

\vspace{2mm}

\textit{Remark.} We develop a non-rigorous qualitative argument that the condition of $(a)$ is nonzero. Note that all terms in $(a)$ can be nonzero (this is only not true in certain circumstances, such as the Christoffel symbols are constant, thus differentiation yields $0$). Thus, we can scale the metric nonlinearly, such as through $g_{ij} \rightarrow F_{ij}(g_{ij})$, if the equation is not automatically nonzero for choice of metric. In general, for the equation $x+y=0, x,y \neq 0$, we can construct a nonlinear function $f$ such that
\begin{equation}
f_1(x) + f_2(y)  \neq 0,
\end{equation}
where $f_1, f_2$ are transformations that all relate to $f$ in some way, such as through the chain rule. For us, $x$ is $(a)$ and $y$ is $(c)$. We also remark our goal is to show there exists such a metric, not that the result holds for all metrics. A nonlinear scaling changes the intrinsic geometry in general.

\vspace{2mm}

\textit{Remark.} Observe the quantity $(a)$ has close connections to the Riemannian tensor if we were to differentiate the Christoffel symbols and rearrange the summations.

\vspace{2mm}

\textit{Remark.} A key of this proof is the fact that the Laplacian of the Ricci tensor along the off-diagonal also depends on the diagonal elements. This is observable with the covariant derivative formula on a $(0,2)$ tensor
\begin{equation}
\label{eqn:cov_deriv_formula_2}
\nabla_l A_{ik} = \partial_l A_{ik} - {\Gamma_{il}}^m A_{mk} - {\Gamma_{kl}}^m A_{im} ,
\end{equation}
and since there are two summations over $m$, and the covariant derivative formula on a $(0,3)$ tensor
\begin{align}
\label{eqn:cov_deriv_formula_3}
\nabla_j A_{lik} = \partial_j A_{lik} -{\Gamma_{jl}}^m A_{mik} - {\Gamma_{ji}}^m A_{lmk} - {\Gamma_{jk}}^m A_{lim} .
\end{align}
As we mentioned, we normally have ${\Gamma_{il}}^k = 0$ for diagonal metrics if $i \neq l \neq k$ \cite{diagonal_christoffel_symbols}.

\section{A special case of the $d$-sphere}
\label{sec:sphere}

Ricci flow upon a sphere has simple closed form. Since curvature of the sphere is already uniformized, Ricci flow upon the sphere has no curvature-steady state to evolve to other than itself. Regardless, contractive properties still exist, so the sphere shrinks until extinction in finite time. The radius of the sphere subject to Ricci flow takes the form
\begin{equation}
r(t) = \sqrt{ R^2 - 2(d-1) t },
\end{equation}
for initial radius $R$ \citep{ricciflowsheridan}. Moreover, the exact coordinate chart mapping from local coordinates to the sphere is known, and governed by the equations
\begin{align}
\mathcal{E}_{\mathbb{S}^{d}}^1 = r(t) \cos u^1,
\mathcal{E}_{\mathbb{S}^{d}}^i = r(t)  \cos u^i  \prod_{j=1}^{i-1} \sin u^j, 
\mathcal{E}_{\mathbb{S}^{d}}^{d+1} = r(t) \prod_{j=1}^{d} \sin u^j .
\end{align}
It is very important to note here that the manifold is embedded in $\mathbb{R}^{d+1}$ dimensional Euclidean space. Here, the local coordinates that parameterize the sphere $u$ are still learned with a neural network. However, we will not work with the coordinate chart. As is observable, it requires high order sinusoidal products. These values are diminished for $\sin u^j < 1$ to some level. We will work directly with a point encoded in space and its normalization. We will take
\begin{align}
\mathcal{E}_{\mathbb{S}^d} = r(t) \frac{u}{||u||_2}
\end{align}
directly.

\vspace{2mm}

Our objective function in this setting is greatly simplified, and the geometric flow loss is
\begin{equation}
\mathcal{L} = \EX  \Big[ \Big|\Big| \mathcal{D}_{\theta_{\mathcal{D}}} ( \mathcal{E}_{\mathbb{S}^{d}} ( u_{\theta_u}(\tilde{\phi}_0), t)) - \tilde{\phi}_{t} \Big|\Big|_2^2 \Big].
\end{equation}
We make some additional remarks. Now, the intrinsic nature of Ricci flow is reduced to an extrinsic parameterization. Rotational vector fields tangent to the sphere, and the translocation of the center of mass no longer exist.

\vspace{2mm}

We can reintroduce some expressivity into the parameterization by introducing a new neural network $\mathcal{S}_{\theta_{\mathcal{S}}} : [0,T] \times \Theta_{\mathcal{S}} \rightarrow \mathbb{R}^{d+1}$. We have
\begin{equation}
\tilde{\mathcal{E}}_{\mathbb{S}^{d-1}}^i = \mathcal{E}_{\mathbb{S}^{d-1}}^i + \mathcal{S}_{\theta_{\mathcal{S}}}^i ,
\end{equation}
where $\tilde{\mathcal{E}}$ is the new sphere after the change of centers. Observe $\mathcal{E}_{\mathbb{S}^{d-1}}$ is the fixed sphere at the center and $\mathcal{S}_{\theta_{\mathcal{S}}}$ is learned in training. We denote $i$ as the $i$-th coordinate in $\mathbb{R}^{d+1}$. In particular, we can decompose the vector field and each particle that composes the manifold evolves according to
\begin{align}
\begin{cases}
\frac{dx(t)}{dt} = \xi_1(x(t)) + \xi_2(x(t)) + \xi_3(x(t))
\\
x(t) = x(0) + \int_0^t ( \xi_1 + \xi_2 + \xi_3 )(x(\tau)) d\tau 
\\
(\xi_1 + \xi_2 + \xi_3 ) \in \Gamma(T\mathcal{M}) \oplus \Gamma(N \mathcal{M})
\end{cases}
,
\end{align}
which correspond to the rotation, translocation, and movement from the geometric flow. Since we cannot easily create rotation, the neural network $\mathcal{S}_{\theta_{\mathcal{S}}}$ means each particle evolves to
\begin{align}
x(t) = x(0) + \int_0^t ( \xi_2 + \xi_3 ) d\tau .
\end{align}

\section{Path integration with Gaussian curvature}
\label{app:gauss_curv_sec}

In this section, we utilize the relation
\begin{equation}
\text{Ric} = K g
\end{equation}
where $K$ is Gaussian curvature, defined as \cite{weisstein_gaussian}
\begin{equation}
K = \frac{1}{\sqrt{\text{det}(g)}} \Bigg( \frac{\partial}{\partial u^2} ( \frac{\sqrt{\text{det}(g)}}{g_{11}} {\Gamma_{11}}^2 ) - \frac{\partial}{\partial u^1}  ( \frac{\sqrt{\text{det}(g)}}{g_{11}} {\Gamma_{12}}^2 ) \Bigg) .
\end{equation}
Thus, our loss function in a learning setting starts out as
\begin{equation}
\EX [ || (\partial_t -  K)g ||_F^2 ] .
\end{equation}
Note that this is exactly Ricci flow, since $\text{Ric}= K g$ for 2-manifolds \cite{mathstackriccigaussiancurv}.

\vspace{2mm}

We wish to apply Green's theorem to simplify $K$ but there is a major issue because $1/\sqrt{\text{det}(g)}$ scales $K$ and we cannot immediately apply Green's theorem in this form. However, we will modify our loss function to allow this. Now, we have 
\allowdisplaybreaks
\begin{align}
& \int \EX \Big[  \sqrt{\text{det}(g)} \Big|\Big| \partial_t g - K g \Big|\Big|_F^2 \Big] dt
\\
= & \int \EX \Big[  \sqrt{\text{det}(g)}  \sum_{ij} (  \partial_t g_{ij} - K g_{ij} )^2 \Big] dt 
\\
= & \int \EX \Big[  \sum_{ij} ( \sqrt[4]{\text{det}(g)} \partial_t g_{ij} - \sqrt[4]{\text{det}(g)} K g_{ij} )^2 \Big]  dt
\\
= & \iiint_{\Sigma \times [0,T]} \sum_{ij} ( \sqrt[4]{\text{det}(g)} \partial_t g_{ij} - \sqrt[4]{\text{det}(g)} K g_{ij} )^2 \rho(u) du^1 du^2 dt
\\
= & \iiint_{\Sigma \times [0,T]} \sum_{ij} ( \sqrt[4]{\text{det}(g)} \partial_t g_{ij} - \sqrt[4]{\text{det}(g)} K g_{ij} )^2 \rho(u) du^1 du^2 dt
\\
= & \iiint_{\Sigma \times [0,T]} \sum_{ij} ( \sqrt{\text{det}(g)} (\partial_t g_{ij})^2 - 2 \sqrt{\text{det}(g)} K g_{ij} \partial_t g_{ij} +  \sqrt{\text{det}(g)} K^2 g_{ij}^2 ) \rho(u) du^1 du^2 dt .
\end{align}
Now, the quadratic Gaussian curvature is problematic because this integral is nontrivial to evaluate and simplify. We propose a first-order Taylor surrogate
\begin{align}
K^2 \approx \widetilde{K}^2 + 2 \widetilde{K}(K - \widetilde{K}) = 2 \widetilde{K} K - \widetilde{K}^2.
\end{align}
Here, $\widetilde{K}$ is a Gaussian curvature approximation that is linear in true Gaussian curvature. Now, Gaussian curvature is approximated linearly. Returning to our loss,
\begin{align}
& \iiint_{\Sigma \times [0,T]} \sum_{ij} ( \sqrt{\text{det}(g)} (\partial_t g_{ij})^2 - 2 \sqrt{\text{det}(g)} K g_{ij} \partial_t g_{ij} +  \sqrt{\text{det}(g)} K^2 g_{ij}^2 ) \rho(u) du^1 du^2 dt 
\\
= & \iiint_{\Sigma \times [0,T]} \sum_{ij} ( \sqrt{\text{det}(g)} (\partial_t g_{ij})^2 - 2 \sqrt{\text{det}(g)} K g_{ij} \partial_t g_{ij} +  \sqrt{\text{det}(g)}  [ \widetilde{K}^2 + 2 \widetilde{K}(K - \widetilde{K}) ]  g_{ij}^2 ) \rho(u) du^1 du^2 dt 
\\
= & \iiint_{\Sigma \times [0,T]} \sum_{ij} ( \sqrt{\text{det}(g)} (\underbrace{ \partial_t g_{ij})^2 }_{\text{(I)} \geq 0} - 2 \sqrt{\text{det}(g)} \underbrace{ ( g_{ij} \partial_t g_{ij} - 2 \widetilde{K} g_{ij}^2 ) K }_{\text{(II)}} -  \sqrt{\text{det}(g)} \underbrace{ \widetilde{K}^2 g_{ij}^2 }_{\text{(III)} \geq 0}   ) \rho(u) du^1 du^2 dt  = 0.
\end{align}
It is of consideration to introduce the change of variables
\begin{align}
d\psi_{ij} = ( g_{ij} \partial_t g_{ij} - 2 \widetilde{K} g_{ij}^2 ) du^1 \wedge du^2 .
\end{align}
This leads to highly rigorous computation, and this route is not preferred. Instead, notice that not only is the integral zero but the integrand should be zero almost everywhere. Thus, it is valid to introduce
\begin{align}
& \iiint_{\Sigma \times [0,T]} \sum_{ij} ( \sqrt{\text{det}(g)} ( \partial_t g_{ij})^2  - 2 \sqrt{\text{det}(g)}  ( g_{ij} \partial_t g_{ij} - 2 \widetilde{K} g_{ij}^2 ) K -  \sqrt{\text{det}(g)} \widetilde{K}^2 g_{ij}^2    ) \rho(u) du^1 du^2 dt 
\\
&  = \iiint_{\Sigma \times [0,T]} \sum_{ij} ( \sqrt{\text{det}(g)} \frac{ ( \partial_t g_{ij})^2 }{( g_{ij} \partial_t g_{ij} - 2 \widetilde{K} g_{ij}^2 )}  - 2 \sqrt{\text{det}(g)}   K -  \sqrt{\text{det}(g)} \frac{ \widetilde{K}^2 g_{ij}^2 }{( g_{ij} \partial_t g_{ij} - 2 \widetilde{K} g_{ij}^2 )}  ) \rho(u) du^1 du^2 dt  .
\end{align}
We have clearly assumed the denominator is nonzero. This is mathematically valid since
\begin{align}
& \EX \Bigg[ \EX \Bigg[   \sum_{ij} \Big( \sqrt{\text{det}(g)} ( \partial_t g_{ij})^2   - 2 \sqrt{\text{det}(g)}  ( g_{ij} \partial_t g_{ij} - 2 \widetilde{K} g_{ij}^2 ) K  -  \sqrt{\text{det}(g)}  \widetilde{K}^2 g_{ij}^2    \Big) \Bigg] \Bigg] = 0
\\
& \implies  \Big( \sqrt{\text{det}(g)} ( \partial_t g_{ij})^2   - 2 \sqrt{\text{det}(g)}  ( g_{ij} \partial_t g_{ij} - 2 \widetilde{K} g_{ij}^2 ) K  -  \sqrt{\text{det}(g)}  \widetilde{K}^2 g_{ij}^2    \Big)   = 0 \ \ \ \text{a.e.}
\end{align}
which follows by nonnegativity, thus division by a nonzero function is perfectly fine. As a side remark, since we used a Gaussian curvature proxy, this introduces error, so the nonnegativity is not guaranteed but holds with a sufficient error tolerance. Now the integration of Gaussian curvature is computable.

\vspace{2mm}

Thus, our Gaussian curvature with the product with the volume element is
\begin{align}
\tilde{K} = \sqrt{\text{det}(g)} K =    \frac{\partial}{\partial u^2} ( \frac{\sqrt{\text{det}(g)}}{g_{11}} {\Gamma_{11}}^2 ) -  \frac{\partial}{\partial u^1}  ( \frac{\sqrt{\text{det}(g)}}{g_{11}} {\Gamma_{12}}^2 )  .
\end{align}
Now, we will use Stokes' / Green's theorem. Observe
\begin{align}
& \iint_{B_{r}(x_0)}  \Bigg( \frac{\partial}{\partial u^2} ( \frac{\sqrt{\text{det}(g)}}{g_{11}} {\Gamma_{11}}^2 ) - \frac{\partial}{\partial u^1}  ( \frac{\sqrt{\text{det}(g)}}{g_{11}} {\Gamma_{12}}^2 ) \Bigg) du^1 du^2
\\
& = \oint_{\partial B_r(x_0)} \frac{\sqrt{\text{det}(g)}}{g_{11}} {\Gamma_{11}}^2 du^1 + \frac{\sqrt{\text{det}(g)}}{g_{11}} {\Gamma_{12}}^2 du^2 
= \oint_{\partial B_r(x_0)} \frac{\sqrt{\text{det}(g)}}{g_{11}} \Bigg( {\Gamma_{11}}^2 du^1 +  {\Gamma_{12}}^2 du^2 \Bigg) .
\end{align}
This allows us to bypass differentiation of the Christoffel symbols, increasing efficiency. Indeed, it is not hard to evaluate this line integral over arbitrary circles, and it is first order in its derivatives with the Christoffel symbols of the second kind. To compute this, we parameterize any circle
\begin{align}
v^1(\omega) = u_0^1 + r \cos \omega, v^2(\omega) = u_0^2 + r \sin \omega 
\end{align}
and make it discrete. We get
\begin{align}
& \oint_{\partial B_r(x_0)} \frac{\sqrt{\text{det}(g)}}{g_{11}} \Bigg( {\Gamma_{11}}^2 du^1 +  {\Gamma_{12}}^2 du^2 \Bigg)
\approx \underbrace{ \sum_{\gamma_i}  \frac{\sqrt{\text{det}(g(v(\omega_i),t))}}{g_{11}} \Big( r {\widetilde{\Gamma}_{11}}^2 ( - \sin \omega_i ) + r {\widetilde{\Gamma}_{12}}^2  ( \cos \omega_i ) \Big) \Delta \omega }_{\text{closed line path loss}} .
\end{align}
Our full discrete loss function is 
\begin{align}
& \mathcal{L} =  \text{closed line path loss}
\\
\label{eqn:loss_gausscurv_remaining_geo_flow}
& + \sum_{(\tilde{\phi}_0^{(i)}, t_i)} \sum_{jk}  \frac{ \sqrt{\text{det}(g(u(\tilde{\phi}_0^{(i)}), t_i))} \Big(  \partial_t g_{jk}^2 (u(\tilde{\phi}_0^{(i)}), t_i) - g_{jk}^2(u(\tilde{\phi}_0^{(i)}), t_i) \widetilde{K}^2 \Big) }{ \Big( g_{ij}(u(\tilde{\phi}_0),t_i) \partial_t g_{ij}(u(\tilde{\phi}_0),t_i) - 2 \widetilde{K} g_{ij}^2(u(\tilde{\phi}_0),t_i) \Big) }
\end{align}
We give some final remarks. We will assign our Gaussian curvature proxy as
\begin{align}
\widetilde{K}_i \leftarrow   \frac{1}{\text{Length}(\partial B_r(x_i))} \oint_{\partial B_r(x_0)} \frac{\sqrt{\text{det}(g)}}{g_{11}} \Bigg( {\Gamma_{11}}^2 du^1 +  {\Gamma_{12}}^2 du^2 \Bigg)   . 
\end{align}
This can be evaluated with a Riemann sum using polar and the increment $r \Delta \theta$. Note that the increment and length denominator have high cancellation. More specifically, $\widetilde{K}$ can be computed as, over a circle,
\begin{align}
\widetilde{K}_i & = \frac{1}{2 \pi r} \sum_i \frac{\sqrt{\text{det}(g_i)}}{g_{11,i}} ( {\widetilde{\Gamma}_{11}}^2 (-\sin \omega_i ) + {\widetilde{\Gamma}_{12}}^2 (\cos \omega_i ) ) r \Delta \omega_i  =  \frac{1}{L} \sum_i \frac{\sqrt{\text{det}(g_i)}}{g_{11,i}} ( {\widetilde{\Gamma}_{11}}^2 (-\sin \omega_i ) + {\widetilde{\Gamma}_{12}}^2 (\cos \omega_i ) )  ,
\end{align}
where $L$ is the number of increments in which the circle of domain $2 \pi$ is divided.

\vspace{2mm}

\textit{Remark.} We remark that this method of simplifying the computation of Gaussian curvature with Stokes' theorem is closely related to the Gauss-Bonnet theorem, which states
\begin{align}
\iint_{\mathcal{M}} K dA = 2 \pi \chi(\mathcal{M}) - \oint_{\partial \mathcal{M}} k_g ds .
\end{align}
Here, $k_g$ is geodesic curvature; however, our method is generally preferred because we do not need to compute a normal vector, which is typical to compute geodesic curvature.

\vspace{2mm}

We make some other remarks. Algebraically, our mathematics is fine, and division by a nonzero function is acceptable. It is important to distinguish that machine learning attempts to minimize, and the true output of the loss is not zero. Thus, the loss landscape is quite different, and the solution behavior across training will be nonuniform among these divisory scalings. Nonetheless, the same solution works if or if not the nonzero function has been divided, and for our purposes, the function division, while affective of the loss landscape, is unimportant for an overall final outcome.

\vspace{2mm}

As a last remark, note that our Gaussian curvature loss requires sampling $\epsilon$-balls. We rely on an integral equivalence
\begin{align}
\iiint \sqrt{\text{det}(g)} K dA dt \propto \int \Bigg( \frac{1}{N} \sum_{i=1}^N \oint_{\partial B_r(x_i)} \frac{\sqrt{\text{det}(g)}}{g_{11}} {\Gamma_{11}}^2 du^1 + \frac{\sqrt{\text{det}(g)}}{g_{11}} {\Gamma_{12}}^2 du^2 \Bigg) dt ,
\end{align}
which is valid under a (non-rigorous) Monte Carlo-type argument. In particular, up to an orientation on the circulation
\begin{align}
& \EX_{x_i \sim U(\Sigma)} [ \sqrt{\text{det}(g)} K |_{x_i} ] \stackrel{\infty \leftarrow N}{\longleftarrow} \frac{1}{N 2 \pi r} \sum_{i=1}^N \oint_{\partial B_r(x_i)} \frac{\sqrt{\text{det}(g)}}{g_{11}} {\Gamma_{11}}^2 du^1 + \frac{\sqrt{\text{det}(g)}}{g_{11}} {\Gamma_{12}}^2 du^2
\\
& \approx \frac{1}{N| \Sigma |} \sum_{i \in  \{ i \in \mathbb{N} : u_i \sim U(\Sigma) \}}^N   | \Sigma | \sqrt{\text{det}(g_i)} K_i \approx \frac{1}{| \Sigma |} \iint_{\Sigma} \sqrt{\text{det}(g)} K dA = \dashint_{\Sigma} \sqrt{\text{det}(g)} K dA .
\end{align}
$| \Sigma |$ denotes Lebesgue measure of $\Sigma$. The first  equation is the expected Gaussian curvature across each ball, the second is the empirical average Gaussian curvature (times volume element), the third is an averaged Monte Carlo Gaussian curvature, and the fourth and fifth are true integrals that compute the average. As a few side remarks: we have omitted weighting function $\rho$ that we saw earlier for simplicity; the Lebesgue measure in the third term indeed cancels, but we include both to illustrate that the following integral is divided by Lebesgue measure (certainly it is nonzero).

\section{New flows with Taylor expansions and curvature-derived tensors}
\label{app:second_FF_sec}

In this section, we develop a custom geometric flow that is based off of Ricci flow, but as we will see, it is also quite different. Ricci flow is widely studied in literature for a variety of theoretical purposes. While some of these purposes apply to a machine learning context, many properties do not. As mentioned in the introduction, a crucial key for us is either contractive or expansive properties. We will derive a computationally-efficient flow and show it satisfies such a property. This flow we will derive has critical extrinsic definitions; however, the PDE will be associated on the Riemannian metric, and no use of vector fields or ordinary differential equations is required.

\vspace{2mm}

The primary value of this section is that we do not need invaluable geometric analytical properties for our contexts, as long as the flow is sufficiently smooth and nice and maintains the properties we desire. Nonetheless, the flow we will create is curvature-dependent, thus falls into the critical class of flows that motivate our work. Indeed, this flow is slightly ad hoc, and that is somewhat the point. We conjecture properties such as hyperbolicity and curvature uniformization, maybe beneficial indirectly, are not key in latent spaces. Our flow is motivated geometrically, while unconventional, and it still achieves good empirical results.

\vspace{2mm}

The all lower-indexed (0,4) Riemannian tensor $R_{ijkl}$ can be decomposed \cite{hamilton1982three}
\begin{align}
R_{ijkl} & = W_{ijkl} + S_{ijkl} + U_{ijkl}
\\
& = W_{ijkl} + \frac{R}{d(d-1)} (g_{ik} g_{jl} -  g_{il} g_{jk} ) +  \frac{1}{d-2}( g_{ik} \widehat{R}_{jl} - g_{il} \widehat{R}_{jk} - g_{jk} \widehat{R}_{il} + g_{jl} \widehat{R}_{ik} ) .
\end{align}
Here, $\widehat{R}_{ij} = \text{Ric}_{ij} - \frac{1}{d} R g_{ij}$ is the traceless part of the Ricci tensor. We will omit the Weyl tensor from our computation since it is trace-free, and we will develop a new tensor that is derivative of
\begin{align}
\widetilde{R}_{ijkl} & = S_{ijkl} + U_{ijkl}
\\
& =   \frac{R}{d(d-1)} (g_{ik} g_{jl} -  g_{il} g_{jk} ) +  \frac{1}{d-2}( g_{ik} \widehat{R}_{jl} - g_{il} \widehat{R}_{jk} - g_{jk} \widehat{R}_{il} + g_{jl} \widehat{R}_{ik} ) .
\end{align}
We will replace the traceless Ricci tensor, which is curvature-based, with a new curvature-type tensor derived from our second fundamental form proxy. From the Gauss equation, it is known
\begin{align}
R_{ijkl} = \overline{R}_{ijkl} + \Pi_{ik} \Pi_{jl} - \Pi_{jk} \Pi_{il} ,
\end{align}
where $h$ is the second fundamental form, which is defined by
\begin{equation}
h(X,Y) = \overline{g} ( \overline{\nabla}_X \nu, Y ) ,
\end{equation}
for tangent vectors $X,Y$, ambient metric $\overline{g}$, Levi-Civita connection $\overline{\nabla}$, and unit normal $\nu$. In local coordinates, we get
\begin{align}
\Pi_{ij} = \overline{g}_{ab} \nu^a \Big( \frac{\partial^2 X^b}{\partial u^i \partial u^j} + {\overline{\Gamma}_{cd}}^b \frac{\partial X^c}{\partial u^i} \frac{\partial X^d}{\partial u^j} \Big) .
\end{align}
Clearly, the Gauss equation is quadratic in second fundamental form, but it is dependent on the second fundamental form nonetheless. Ricci curvature follows by taking the trace, and a linear version is reasonable under Taylor expansions and in special cases. As a general rule,
\begin{align}
f(\Pi) = f(\Pi_0) + f'(\Pi_0) \delta \Pi + O(\delta \Pi^2)
\end{align}
where we have Taylor expanded $f$. Returning to the Gauss equation, and denoting $f$ as the quadratic difference in second fundamental form, we get with a Taylor expansion
\begin{align}
f(\Pi) = f(\Pi_0) + \underbrace{ \Bigg( \Pi_{ik,0} \delta \Pi_{jl} + \delta \Pi_{ik} \Pi_{jl,0} - \Pi_{jk,0} \delta \Pi_{il} - \delta \Pi_{jk} \Pi_{il,0} \Bigg) }_{= \text{ linear in } \Pi_0, \delta \Pi} + \text{higher order terms} .
\end{align}
Here, we have used the fact that
\begin{align}
\sum_{ij} \frac{ \partial f }{\partial \Pi_{ij}} \Bigg|_{\Pi_0}  \delta \Pi_{ij} = \Pi_{ik,0} \delta \Pi_{jl} + \delta \Pi_{ik} \Pi_{jl,0} - \Pi_{jk,0} \delta \Pi_{il} - \delta \Pi_{jk} \Pi_{il,0} ,
\end{align}
since each term in each quadratic is indexed differently. Note that the corresponding tensor differentiation rule is found in \cite{Itskov}. Thus, the Gauss equation holds as a linear equation in certain circumstances. Here, $\Pi_0$ is the value in which the series is expanded. As an additional remark, the only other higher order term is quadratic, since it is a Taylor expansion of a quadratic function, i.e. $\partial^3 f$ vanishes.

\vspace{2mm}

In our experiments, we have initialized our geometry at $t=0$ as the sphere, which is constant in curvature. Hence, the condition
\begin{align}
\Pi_{ij,0} = \kappa g_{ij,0}
\end{align}
is reasonable for curvature $\kappa$. Thus, our Taylor expanded linear term around a perturbation of the sphere at initial time gives
\begin{align}
\Pi_{ik,0} \delta \Pi_{jl} + \delta \Pi_{ik} \Pi_{jl,0} - \Pi_{jk,0} \delta \Pi_{il} - \delta \Pi_{jk} \Pi_{il,0}   =  \kappa \Bigg( g_{ik,0} \delta \Pi_{jl} + \delta \Pi_{ik} g_{jl,0} - g_{jk,0} \delta \Pi_{il} - \delta \Pi_{jk} g_{il,0} \Bigg)  .
\end{align}
Contracting,
\begin{align}
\kappa g^{jl,0} \Bigg( g_{ik,0} \delta \Pi_{jl} + \delta \Pi_{ik} g_{jl,0} - g_{jk,0} \delta \Pi_{il} - \delta \Pi_{jk} g_{il,0} \Bigg) = \kappa \Bigg(    g_{ik,0} \text{Tr} ( \delta \Pi ) + d \delta  \Pi_{ik} -\delta \Pi_{ik} - \delta \Pi_{ik} \Bigg)  .
\end{align}
Adding the above back to our quadratic difference and contracting,
\begin{align}
\text{Ric}_{ik} & = g^{jl,0} \Bigg( \Pi_{ik,0} \Pi_{jl,0} - \Pi_{jk,0} \Pi_{il,0} \Bigg) + \kappa g_{ik,0} \text{Tr}(\delta \Pi) + \kappa(d-2)\delta \Pi_{ik} + \mathcal{O}(\delta \Pi^2) \\
& =  \kappa^2 (d-1) g_{ik,0}  + \kappa g_{ik,0} \text{Tr}(\delta \Pi) + \kappa(d-2)\delta \Pi_{ik} + \mathcal{O}(\delta \Pi^2) 
 \\
& =  \kappa (d-1) \Pi_{ik,0}  + \kappa g_{ik,0} \text{Tr}(\delta \Pi) + \kappa(d-2)\delta \Pi_{ik} + \mathcal{O}(\delta \Pi^2)
 \\
& =  \kappa \Bigg[ (d-1) \Pi_{ik,0}  +  g_{ik,0} \text{Tr}(\delta \Pi) + (d-2)\delta \Pi_{ik} \Bigg] + \mathcal{O}(\delta \Pi^2).
\end{align}
Dropping the high-order terms, we see sufficiently close to the sphere
\begin{align}
\delta \text{Ric}_{ik} \approx \kappa g_{ik,0} \text{Tr}(\delta \Pi) + \kappa(d-2)\delta \Pi_{ik} .
\end{align}
Therefore, the perturbation in Ricci curvature is proportional up to a perturbation and constants in $\Pi$ in our framework, as desired. Moreover, $\text{Ric}$ is proportional to $\Pi$ at $t=0$ as well.

\vspace{2mm}

The second fundamental form is expensive to compute, so we will use a second fundamental form surrogate
\begin{align}
\Pi_{ij} = \Bigg| \Bigg| \partial_{ij} \mathcal{E} - \sum_{kl} g^{kl} \langle \partial_{ij} \mathcal{E}, \partial_k \mathcal{E} \rangle \partial_l \mathcal{E} \Bigg| \Bigg|_2 .
\end{align}
We will specifically restrict our manifolds embedded in Euclidean space, which implies that $\overline{R}_{ijkl} = 0$, and so we get the simplification
\begin{align}
\Pi_{ij} = \partial_{ij} X \cdot \nu .
\end{align}
Here, $X$ is the parametric map mapping from local coordinates to the embedded manifold. Thus, our motivation for the design of this new tensor is based on the Gauss equation. In particular, the second fundamental form has close connections to Ricci curvature, thus we conciliate this notions to define the following tensor 
\begin{align}
\widetilde{H}_{ijkl} = \frac{H}{d(d-1)} (g_{ik} g_{jl} -  g_{il} g_{jk} ) +  \frac{1}{d-2}( g_{ik} \widehat{H}_{jl} - g_{il} \widehat{H}_{jk} - g_{jk} \widehat{H}_{il} + g_{jl} \widehat{H}_{ik} ) ,
\end{align}
where we will set
\begin{align}
\widehat{H}_{ij} = \Pi_{ij} - \frac{H}{d} g_{ij} ,
\end{align}
where $H$ is mean curvature. The exchange of $R$ and $H$ is justified since
\begin{align}
R & = g^{ij} \text{Ric}_{ij} = g^{ij} \kappa^2 (d-1) g_{ij}  + \mathcal{O}(\delta \Pi) = \kappa^2 d (d-1)  + \mathcal{O}(\delta \Pi) 
\end{align}
and since
\begin{align}
H = g^{ij} \Pi_{ij} = g^{ij} (\kappa g_{ij}) = \kappa d .
\end{align}
But recall in our Taylor expansion derivation $\text{Ric}$ was quadratic in $\kappa$ and $\Pi$ is linear in $\kappa$, but $\delta \text{Ric}$ is linear in $\kappa$ with respect to $\delta \Pi$. Thus, the missing factor of $\kappa$ justifies the replacement. It can be noted both
\begin{align}
\widehat{R}_{ij} = \text{Ric}_{ij} - \frac{1}{d} R g_{ij} = \widehat{H}_{ij} = \Pi_{ij} - \frac{H}{d} g_{ij}  = 0 
\end{align}
in the case of the sphere, and so ours (approximated) is consistent with its derivation (true).

\vspace{2mm}

One may consider an addition of a Weyl tensor replacement. In our experiments, we consider this. Our final tensor is
\begin{align}
H_{ijkl} =  \frac{H}{d(d-1)} (g_{ik} g_{jl} -  g_{il} g_{jk} )  +  \frac{1}{d-2}( g_{ik} \widehat{H}_{jl} - g_{il} \widehat{H}_{jk} - g_{jk} \widehat{H}_{il} + g_{jl} \widehat{H}_{ik} )   .
\end{align}
From this tensor, we can build the geometric flow
\begin{align}
\label{eqn:custom_geo_flow}
\partial_t g_{ij} = H_{ij} = g^{kl} H_{ikjl} .
\end{align}

\vspace{2mm}

\textbf{Theorem 2.} The geometric flow of equation \ref{eqn:custom_geo_flow} is expansive or contractive in its volume element, i.e. $\partial_t \sqrt{\text{det}(g)} \neq 0$, if \begin{align}
 H    \neq 0  ,    
\end{align}
assuming it is nonsingular.

\vspace{2mm}

\textit{Proof.}  Observe
\begin{align}
& \partial_t \sqrt{\text{det}(g)} = \frac{1}{2} \sqrt{\text{det}(g)} g^{ij} \partial_t g_{ij} 
\\
& = \frac{1}{2} \sqrt{\text{det}(g)} g^{ij}  H_{ij} 
\\
& = \frac{1}{2} \sqrt{\text{det}(g)} g^{ij}  g^{kl}  \Bigg( \frac{H}{d(d-1)} (g_{ij} g_{kl} -  g_{il} g_{kj} ) +  \frac{1}{d-2}( g_{ij} \widehat{H}_{kl} - g_{il} \widehat{H}_{kj} - g_{jk} \widehat{H}_{il} + g_{kl} \widehat{H}_{ij} )  \Bigg)
\\
& = \frac{1}{2} \sqrt{\text{det}(g)} \Bigg(  \frac{H}{d(d-1)} (d^2 -  \delta_{i}^k \delta_k^i ) +  \frac{1}{d-2}( d g^{kl} \widehat{H}_{kl} - \delta_i^k g^{ij} \widehat{H}_{kj} - \delta_k^i g^{kl} \widehat{H}_{il} + d g^{ij} \widehat{H}_{ij} )  \Bigg)
\\
& = \frac{1}{2} \sqrt{\text{det}(g)} \Bigg(  \frac{H}{d(d-1)} (d^2 -  d ) +  \frac{1}{d-2}( d g^{kl} \widehat{H}_{kl} - \text{Tr}_g ( \widehat{H}) - \text{Tr}_g(\widehat{H}) + d g^{ij} \widehat{H}_{ij} )  \Bigg)
\\
& = \frac{1}{2}\sqrt{\text{det}(g)} \Bigg(  \frac{H}{d(d-1)} (d^2 -  d ) +  \frac{1}{d-2}( 2(d-1) \text{Tr}_g(\widehat{H}) )  \Bigg)
\\
& = \frac{1}{2} \sqrt{\text{det}(g)} \Bigg( \frac{H}{d(d-1)} d(d -  1 )  \Bigg) =   \frac{1}{2} \sqrt{\text{det}(g)} (  H  ) .
\end{align}
The first line is by Jacobi's formula.

$ \square $

\vspace{2mm}

\textbf{Theorem 3.} Suppose $g \in \Gamma(\text{Sym}^2 T^* \mathcal{M})$ for a manifold of intrinsic dimension $3$ or greater, and that $||g||_F \rightarrow 0$ with respect to some finite time. Assume a regularity condition on the decay rate of $g$. Then there exists a Riemannian metric $g$ such that $|g^{kl} H_{ikjl}| \not\to 0$.

\vspace{2mm}

\textit{Proof.} Notice
\begin{align}
& | g^{kl} H_{ikjl} |  = \Bigg| g^{kl}  \Bigg(  \frac{H}{d(d-1)} (g_{ij} g_{kl} -  g_{il} g_{kj} ) +  \frac{1}{d-2}( g_{ij} \widehat{H}_{kl} - g_{il} \widehat{H}_{kj} - g_{jk} \widehat{H}_{il} + g_{kl} \widehat{H}_{ij} ) \Bigg) \Bigg|
\\
& = \Bigg|  \frac{H}{d(d-1)} ( dg_{ij}  -  g_{il} \delta_j^l ) +  \frac{1}{d-2}( g_{ij} g^{kl} ( \Pi_{kl} - \frac{H}{d} g_{kl} ) - \delta_i^k \widehat{H}_{kj} - \delta_j^l \widehat{H}_{il} + d \widehat{H}_{ij} )   \Bigg|
\\
& = \Bigg|  \frac{H}{d(d-1)} ( dg_{ij}  -  g_{ij} ) +  \frac{1}{d-2}( g_{ij} g^{kl}  \Pi_{kl} - H g_{ij}   - 2 \widehat{H}_{ij}  + d \widehat{H}_{ij} )   \Bigg| .
\end{align}
Now taking the limit, and using an equivalence $||g||\rightarrow 0 \implies |g_{ij}| \rightarrow 0$ under regularity,
\begin{align}
& \lim_{|g_{ij}|\rightarrow 0 \ \forall i,j} \Bigg|  \frac{H}{d(d-1)} ( dg_{ij}  -  g_{ij} ) +  \frac{1}{d-2}( g_{ij} g^{kl}  \Pi_{kl} - H g_{ij}   - 2 \widehat{H}_{ij}  + d \widehat{H}_{ij} ) \Bigg| 
\\
& = 
\lim_{|g_{ij}|\rightarrow 0 \ \forall i,j} \Bigg|    \frac{1}{d-2}( g_{ij} g^{kl}  \Pi_{kl} - H g_{ij}   - 2 \widehat{H}_{ij}  + d \widehat{H}_{ij} )   \Bigg| = 
\lim_{|g_{ij}|\rightarrow 0 \ \forall i,j} \Bigg| \widehat{H}_{ij}  \Bigg|  = 
\lim_{|g_{ij}|\rightarrow 0 \ \forall i,j} \Bigg| \Pi_{ij} \Bigg| .
\end{align}
In general, this is not identically zero and we have the result.

$ \square $

\vspace{2mm}

We will compute mean curvature as
\begin{align}
H_{\text{approx}} & \propto  \Big| \Big| \Delta_g \mathcal{E} \Big| \Big| = 
\Bigg| \Bigg| \begin{pmatrix}
\Delta_g \mathcal{E}_1 \\ \vdots \\ \Delta_g \mathcal{E}_{D}    
\end{pmatrix}
\Bigg| \Bigg|, \ \ \ \ \ \text{or equivalently,} \ \ \ \ \ \Delta_{\Sigma} \mathcal{E}_i = - H n_i.
\end{align}
We have ignored constants, as sometimes the equation on the right is rescaled under conventions. This will be our approximation to mean curvature. Again, $\mathcal{E}$ is the parametric map, and the Laplace-Beltrami operator is applied element-wise.

\section{Discussion of new non-parametric geometric flows with special cases of Perelman's $\mathcal{W}$-functional} 
\label{app:perelman_functional}

We discuss use of functionals for our methods with Perelman's functional. First, we make the following remark about physics-informed neural networks (PINNs). We can simulate a physics loss by minimizing a PDE residual over a manifold which has a corresponding integral formulation of the form
\begin{equation}
|| \ f \ | |_{L^1(\mathcal{M}_t \times [0,T])} \propto \int_0^T \int_{\mathcal{M}_t} |f| \tilde{\rho} dV_t dt = \int_0^T \int_{\Sigma} |  f \circ \mathcal{E}(u,t) | \rho(u) \sqrt{\text{det}(g)} \bigwedge_k du^kdt 
\end{equation}
for suitable Radon-Nikodym derivative $\rho$ with corresponding weighting function $\tilde{\rho}$ over the manifold, and the physics loss can be evaluated discretely with a collocation procedure. In a PINN setup, $f$ is a PDE residual, but this need not be the case. Even moreso, we can simulate integration over a manifold with the same sampling procedure as a PINN by utilizing Monte-Carlo integration, i.e.
\begin{equation}
\lim_{N \rightarrow \infty} \frac{ \int_0^T \int_{\Sigma} \sqrt{\text{det}(g(u,t))} \bigwedge_k du^kdt }{N} \sum_{u_i \in U(\rho), 1 \leq i \leq N} |f(\mathcal{E}(u_i,t_i))| \sqrt{\text{det}(g(u_i, t_i))} \overset{P}{\rightarrow} \int_0^T \int_{\mathcal{M}_t}  |f(x)| \tilde{\rho} dV_t dt 
\end{equation}
by the (weak) law of large numbers, and we have smooth manifold map $\mathcal{E} : \Sigma \rightarrow \mathcal{M}_t$ and function over the manifold $f : \mathcal{M}_t \rightarrow \mathbb{R}^d$. In practice, $\int_0^T \int_{\Sigma} \sqrt{\text{det}(g(u,t))} \bigwedge_k du^kdt$ is constant and irrelevant for machine learning purposes. 

\vspace{2mm}

Perelman introduced the $\mathcal{W}$-functional \cite{perelman2002entropyformularicciflow} \cite{toppingricciflow}
\begin{equation}
\label{eqn:perelman_functional}
\mathcal{W}(f,g,\tau) = \int [ \tau ( R + || \nabla f ||^2) + f - n ] u dV
\end{equation}
for function $f$, metric $g$, scale parameter $\tau$, and $u = (4 \pi \tau)^{-n/2} e^{-f}$. Here, $R$ is scalar curvature. Under Ricci flow as well as other conditions, it can be shown Ricci flow is a gradient flow and that $\frac{d}{dt} \mathcal{W} \geq 0$ \cite{toppingricciflow}.

\vspace{2mm}

We propose considering the conformally changed metric
\begin{equation}
g(u,t) =  \psi^{4/(d-2)}(u,t) \cdot g 
\end{equation}
with special case of scalar curvature 
\begin{align}
R(\psi^{4/(d-2)}g) = \frac{ 4 \frac{d-1}{d-2} \Delta_g \psi - R(g) \psi}{ \psi^{ (d+2)/(d-2)}} .
\end{align}
Observe there is an equivalent form
\begin{align}
R & = e^{-2 \varphi} (R_g - 2(d-1) \Delta_g \varphi - (d-2)(d-1)g(d \varphi,d \varphi) )
\end{align}
but this form is more expensive computationally due to the $g(d\varphi,d\varphi)$ term, which requires many derivatives. Here, $\Delta_g$ is the Laplace-Beltrami operator. Thus, we propose setting $f$ as fixed and using this closed form scalar curvature in the functional of equation \ref{eqn:perelman_functional} and simulating the geometric flow, which is not necessarily Ricci flow, in a Monte-Carlo integration type objective. We solve
\begin{equation}
\frac{d}{dt}\mathcal{W}_{\text{modified}}  \geq 0 
\end{equation}
with an objective of the form
\begin{align}
&  \Bigg| \text{positive const} -  \frac{d}{dt} \mathcal{W}_{\text{modified}}(0,g,1) \Bigg| =  \Bigg| 
 \text{positive const} -    \frac{d}{dt} \int_{\mathcal{M}_t} R \tilde{\rho} dV_t  \Bigg|
\\
& =  \Bigg| \text{positive const} -  \frac{d}{dt} \int_{\Sigma}  \frac{ 4 \frac{d-1}{d-2} \Delta_g \psi - R(g) \psi}{ \psi^{ (d+2)/(d-2)}} \sqrt{\text{det}(g)} \rho(u) \bigwedge_k du^k \Bigg|
\end{align}
with a Monte-Carlo integration-type loss, where $\psi = \psi_{\theta_{\psi}}$ is a neural network. Here, we compute
\allowdisplaybreaks
\begin{align}
& \frac{d}{dt} \int_{\Sigma} e^{-2 \varphi} \frac{ 4 \frac{d-1}{d-2} \Delta_g \psi - R(g) \psi}{ \psi^{ (d+2)/(d-2)}} \sqrt{\text{det}(g(u,t))} \rho(u) \bigwedge_k du^k\\
\approx
\ & \frac{d}{dt} \sum_{(u, t_i ), (\tilde{\phi}_0^{(i)}, t_i)  \sim (\Phi,[0,T])} \frac{ 4 \frac{d-1}{d-2} \Delta_g \psi - R(g) \psi}{ \psi^{ (d+2)/(d-2)}} \sqrt{\text{det}(g(u,t))} \\
=
& \sum_{(u, t_i ), (\tilde{\phi}_0^{(i)}, t_i)  \sim (\Phi,[0,T])} \frac{d}{dt} \Bigg( \frac{ 4 \frac{d-1}{d-2} \Delta_g \psi - R(g) \psi}{ \psi^{ (d+2)/(d-2)}} \sqrt{\text{det}(g(u,t))} \Bigg) .
\end{align}
Note that the time derivative can be taken with automatic differentiation, and we assume exchange of differentiation and integration is valid. Although this objective forces the metric to undergo evolution (observe $\psi=\text{constant}$ has loss $\gg 0$), we do not have a closed form for the evolution, thus we describe the corresponding $g$ as non-parametric. We remark we choose $f=0$, but $u$ and $\tilde{\rho}$ actually have some relation, and so we simplify the above. Lastly, we also remark knowing $\rho$ is irrelevant, as this weighting matches uniform sampling with respect to training data $\phi_0 \sim \Phi$, as discussed in section \ref{sec:methods}.

\vspace{2mm}

We argue the computational strengths of this method. Note that our evaluation has the quantity
\begin{align}
\partial_t \Delta_g \psi
\end{align}
as its most rigorous computational. Moreover, we have restricted our metric to be conformally changed. A more classical geometric flow would be that of Ricci flow written by the conformal change of the Ricci tensor
\begin{align}
\partial_t g = - (d-2) (\text{Hess}(\psi) - d\psi \otimes d\psi ) - (\Delta \psi + (d-2)|\nabla \psi|^2 ) .
\end{align}
Here, $\text{Hess}$ is the Riemannian Hessian with elements
\begin{align}
( \text{Hess}(\psi) )_{ij} = \frac{\partial^2 \psi}{\partial u^i \partial u^j} - {\Gamma_{ij}}^k \frac{\partial \psi}{\partial u^k} .
\end{align}
where we have used Einstein notation. $\Delta_{\tilde{g}}$ is the Laplace-Beltrami operator
\begin{align}
\Delta_{g} \psi = \frac{1}{\sqrt{\text{det}(g)}} \frac{\partial}{\partial u^i} \Bigg( \sqrt{\text{det}(g)} g^{ij} \frac{\partial \psi}{\partial u^j} \Bigg)  .
\end{align}
Also, we have used
\begin{align}
|\nabla \psi|^2 = g^{ij} \partial_i \psi \partial_j \psi ,
\end{align}
and $(d\psi \otimes d\psi)_{ij} = \partial_i \psi \partial_j \psi$.
As is observable, this method requires far more differentiation and is overall more costly. Our method is more efficient, while non-parametric.

\vspace{2mm}

we provide a derivation of the product rule on the discretized functional as in section \ref{app:perelman_functional}. This section is mostly routine computation, but we provide the details for clarity. Recall we are interested in the quantity
\begin{align}
& \frac{d}{dt} \Bigg( \frac{ 4 \frac{d-1}{d-2} \Delta_g \psi - R(g) \psi}{ \psi^{ (d+2)/(d-2)}} \sqrt{\text{det}(g(u,t))} \Bigg)
\\[2em]
& = \frac{d}{dt} \Big( 4 \frac{d-1}{d-2} \Delta_g \psi - R(g) \psi \Big) \Big( \psi^{ (d+2)/(d-2)}  \Big)^{-1} \sqrt{\text{det}(g(u,t))} 
\\
& \ \ \ +  \Big(4 \frac{d-1}{d-2} \Delta_g \psi - R(g) \psi \Big)  \frac{d}{dt} \Big( \psi^{ (d+2)/(d-2)} \Big)^{-1} \sqrt{\text{det}(g(u,t))}  
\\
& \ \ \ +  \Big( 4 \frac{d-1}{d-2} \Delta_g \psi - R(g) \psi \Big) \Big( \psi^{ (d+2)/(d-2)}  \Big)^{-1} \frac{d}{dt} \sqrt{\text{det}(g(u,t))} 
\\[2em]
& =  \Big(  4 \frac{d-1}{d-2} \partial_t \Delta_g \psi - R(g) \partial_t \psi  \Big) \Big( \psi^{ (d+2)/(d-2)}  \Big)^{-1} \sqrt{\text{det}(g(u,t))} 
\\
& \ \ \ -  \Big( 4 \frac{d-1}{d-2} \Delta_g \psi - R(g) \psi \Big)  \Big( \psi^{ (d+2)/(d-2)} \Big)^{-2} \Big( \frac{d+2}{d-2} \psi^{ (d+2)/(d-2) - 1} \partial_t \psi \Big)
\sqrt{\text{det}(g(u,t))}  
\\
& \ \ \ +  \Big( 4 \frac{d-1}{d-2} \Delta_g \psi - R(g) \psi \Big) \Big( \psi^{ (d+2)/(d-2)}  \Big)^{-1} \Big(  \sqrt{\text{det}(g(u,t))}  \cdot \text{Tr} \Big( g^{-1}(u,t) \partial_t g(u,t) \Big) \Big) .
\end{align}
The very last term is by Jacobi's formula. As an additional remark, our exchange of differentiation is supported by Clairaut's theorem ($\psi$ is a neural network and it is continuous; we will allow invertibility of the metric, so the Laplace-Beltrami operator does not diverge, since it requires dividing by the volume element).

\vspace{2mm}

\textbf{Theorem 4.} Suppose the loss
\begin{align}
\Bigg| \text{positive constant} -  \frac{d}{dt} \int_{\Sigma}  \frac{ 4 \frac{d-1}{d-2} \Delta_g \psi - R(g_0) \psi}{ \psi^{ (d+2)/(d-2)}} \sqrt{\text{det}(g)} \rho(u) \bigwedge_k du^k \Bigg|
\end{align}
is exactly satisfied, and that $\psi$ is a bounded neural network with $C^{\infty}$ activation on $\Sigma \times [0,T]$.  Suppose $g_{ij} \in C^{\infty}$ with respect to time, and that the intrinsic dimension satisfies $d \geq 3$. Suppose sufficient regularity on the eigenvalues, i.e. they converge at a similar rate in the limit. Then there is a constant $\delta>0$ such that $||g||_F > \delta$ for all $t$, i.e. $||g||_F \rightarrow 0$ is impossible. Moreover, we can show this without using the relation between $\psi$ and $g$.

\vspace{2mm}

\textit{Proof.} Using continuity and piecewise linearity,
\begin{align}
& \lim_{||g||_F \rightarrow 0}  \Bigg| C - \Big(  4 \frac{d-1}{d-2} \partial_t \Delta_g \psi - R(g) \partial_t \psi  \Big) \Big( \psi^{ (d+2)/(d-2)}  \Big)^{-1} \sqrt{\text{det}(g(u,t))} 
\\
& -  \Big( 4 \frac{d-1}{d-2} \Delta_g \psi - R(g) \psi \Big)  \Big( \psi^{ (d+2)/(d-2)} \Big)^{-2} \Big( \frac{d+2}{d-2} \psi^{ (d+2)/(d-2) - 1} \partial_t \psi \Big)
\sqrt{\text{det}(g(u,t))}  
\\
&  +  \Big( 4 \frac{d-1}{d-2} \Delta_g \psi - R(g) \psi \Big) \Big( \psi^{ (d+2)/(d-2)}  \Big)^{-1} \Big(  \sqrt{\text{det}(g(u,t))}  \cdot \text{Tr} \Big( g^{-1}(u,t) \partial_t g(u,t) \Big) \Big)\Bigg|
\\[2em]
& =    \Bigg| C - \Big(  4 \frac{d-1}{d-2} \partial_t \Delta_g \psi - R(g) \partial_t \psi  \Big) \Big( \psi^{ (d+2)/(d-2)}  \Big)^{-1} \lim_{||g||_F \rightarrow 0} \sqrt{\text{det}(g(u,t))} 
\\
& -  \Big( 4 \frac{d-1}{d-2} \Delta_g \psi - R(g) \psi \Big)  \Big( \psi^{ (d+2)/(d-2)} \Big)^{-2} \Big( \frac{d+2}{d-2} \psi^{ (d+2)/(d-2) - 1} \partial_t \psi \Big)
\lim_{||g||_F \rightarrow 0} \sqrt{\text{det}(g(u,t))}  
\\
&  +  \Big( 4 \frac{d-1}{d-2} \Delta_g \psi - R(g) \psi \Big) \Big( \psi^{ (d+2)/(d-2)}  \Big)^{-1} \Big(  \lim_{||g||_F \rightarrow 0} \sqrt{\text{det}(g(u,t))}  \cdot \text{Tr} \Big( g^{-1}(u,t) \partial_t g(u,t) \Big) \Big)\Bigg|
\\[2em]
& = \Bigg| C - \Big( 4 \frac{d-1}{d-2} \Delta_g \psi - R(g) \psi \Big) \Big( \psi^{ (d+2)/(d-2)}  \Big)^{-1} \Big(  \lim_{||g||_F \rightarrow 0} \sqrt{\text{det}(g(u,t))}  \cdot \text{Tr} \Big( g^{-1}(u,t) \partial_t g(u,t) \Big) \Big)\Bigg| .
\end{align}
It remains to show
\begin{align}
\lim_{||g||_F \rightarrow 0}  \Bigg| \sqrt{\text{det}(g(u,t))}  \cdot \text{Tr} \Big( g^{-1}(u,t) \partial_t g(u,t) \Big) \Bigg|  > 0  
\end{align}
is not possible. Observe
\begin{align}
& \lim_{||g||_F \rightarrow 0}  \Bigg| \sqrt{\text{det}(g(u,t))}  \cdot \text{Tr} \Big( g^{-1}(u,t) \partial_t g(u,t) \Big) \Bigg|  
\\ 
& \leq \lim_{||g||_F \rightarrow 0}  \Bigg| \sqrt{ \prod_i \lambda_i }  \cdot || g^{-1}(u,t) ||_F ||\partial_t g||_F \Bigg| 
\\
& \leq \lim_{||g||_F \rightarrow 0}  \Bigg| M \sqrt{ \prod_i \lambda_i }  \cdot \frac{\sqrt{d}}{\lambda_1}  \Bigg| 
\\
& = 0 .
\end{align}
Here, we have assumed a bound on the derivative of $g$ using local smoothness since $\Sigma \times [0,T]$ is compact. The bound on the metric inverse derives from a singular value bound (recall the connections between singular values and eigenvalues; the exact bound is inversely linear in singular value, which is equivalent to the eigenvalue since $g$ is positive semi-definite, hence has nonnegative eigenvalues. In fact, we will assume positive definiteness so it is not ill-defined, but metrics can have zero determinants on sets of measure zero). Lastly, the trace norm is a simple application of Cauchy-Scwharz. Clearly if $d>2$, the product outscales the negative first power eigenvalue and we have the equality with zero. Thus, it is impossible for the loss to be satisfied, which completes the proof that $||g||_F \rightarrow 0$ is impossible if the loss is satisfied.

\vspace{2mm}

In the above, it is crucial to assume regularity on the decay rates of the eigenvalues. In particular, the above is not necessarily true without a regularity condition on the decay rates. In a machine learning context, this is reasonable.

$ \square $

\vspace{2mm}

As a last remark, observe that
\begin{align}
&  \frac{ 4 \frac{d-1}{d-2} \Delta_g \psi - R(g) \psi}{ \psi^{ (d+2)/(d-2)}} \sqrt{\text{det}(g(u,t))}  
\\
& =  \frac{ 4 \frac{d-1}{d-2} \Delta_g \psi - R(g) \psi}{ \psi^{ (d+2)/(d-2)}} \psi^{2d / (d-2)} \sqrt{\text{det}(g_0(u))} 
\\
& =  (4 \frac{d-1}{d-2} \Delta_g \psi - R(g) \psi) \psi  \sqrt{\text{det}(g_0(u))} .
\end{align}
Here we have noted
\begin{align}
\psi^{- (d+2)/(d-2)} \cdot \psi^{2d/d-2} =  \psi^{2d/(d-2) - (d+2)/(d-2)}  = \psi .
\end{align}
Note that in our proof of the Theorem, we do not need to assume $\psi \rightarrow 0$, only the restriction on $g$, but they are equivalent up to assumptions on $g$.

\section{Energy functionals with harmonic maps}
\label{app:harmonic_maps}

In this section, we consider geometric flows induced by neural discovery of harmonic maps. A harmonic map stems from a harmonic function between manifolds \cite{HeleinHarmonic}. This section mostly ignores the harmonic aspect of the map, thus this flow is our most heuristic of the four we develop here. As we have discussed, this theoretical effect is of no means to us, thus we can primarily ignore any desire of enforcing our functions to be harmonic.  A harmonic map is a function between manifolds
\begin{align}
\tilde{\psi} : \mathcal{M} \rightarrow \mathcal{N} 
\end{align}
such that the critical point of the energy functional is satisfied. For our purposes, we will restrict manifold $\mathcal{M}= \mathcal{M}_0$ as an initialized geometry with no temporal dependence, and we will take $\mathcal{N} = \mathcal{N}_t$ to be time-evolving, thus our harmonic map is such that
\begin{align}
\tilde{\psi} : \mathcal{M}_0 \rightarrow \mathcal{N}_t .
\end{align}
By abuse of notation, we will denote $\tilde{\psi} = \psi$, as we will define $\psi$ in a moment to be a composition of functions based on the harmonic map with respect to the local coordinate of $\mathcal{M}$. It is known the harmonic map is the critical point in the first variation sense of the energy functional
\begin{align}
E(\psi)(t) = \int_{\Sigma} |d \psi|^2 \bigwedge_k du^k.
\end{align}
In practice, this differential can be computed as
\begin{align}
|d\psi|^2 = g^{ij}(u) h_{ab}(\psi(u)) \partial_i \psi^a \partial_j \psi^b 
\end{align}
for $d\psi \in C^{\infty} (T^* \mathcal{M} \otimes f^* T \mathcal{N}_t)$. Thus, the energy functional in an intrinsic setting
\begin{align}
E(\psi)(t) = \int | d \psi |^2 \star 1 & = \int_{\mathcal{N}_t} g^{ij}(u) h_{ab}(\psi(u)) \frac{\partial \psi^a}{\partial u^i} \frac{\partial \psi^b}{\partial u^j} \text{Vol}_g
\\
& = \int_{\Sigma} g^{ij}(u) h_{ab}(\psi(u)) \frac{\partial \psi^a}{\partial u^i} \frac{\partial \psi^b}{\partial u^j}  \sqrt{\text{det}(g)} \bigwedge_k du^k .
\end{align}
Here, $\star$ is the Hodge star. Notice the above formulation is meaningless in an extrinsic setting, which is necessary for us. For example, the indices $ab$ run over intrinsic dimension of $\mathcal{N}_t$, which does not work for us. We will restrict our attention to the functional when $\psi \in \mathbb{R}^D$ exists in an immersion and characterizes a parametric map. Thus, we take
\begin{align}
\widetilde{E}(\psi)(t) & = \int g^{ij}(u) \Big\langle \frac{\partial \psi}{\partial u^i}, \frac{\partial \psi}{\partial u^j} \Big\rangle_{\mathbb{R}^D} \text{Vol}_g
\\
& = \int_{\Sigma} g^{ij}(u) \Big\langle \frac{\partial \tilde{\psi}(\varphi(u))}{\partial u^i}, \frac{\partial \tilde{\psi}(\varphi(u))}{\partial u^j}  \Big\rangle_{\mathbb{R}^D} \sqrt{\text{det}(g)} \bigwedge_k du^k .
\end{align}
We have set $\varphi$ to be the parametric map from $\Sigma$ to $\mathcal{M}$. We will denote 
\begin{align}
\psi(u) = \tilde{\psi} \circ \varphi(u)    
\end{align}
for short. In a very typical setting, the time derivative of this energy is set to $0$, which is based on the calculus of variations Euler-Lagrange problem. We will not work with the calculus of variations problem because, neurally, it is solved trivially. As we have demonstrated, triviality or near triviality are not ideal for us (see Figure \ref{fig:geo_flow_comparison}). As in \ref{app:perelman_functional}, we will induce a flow by differentiating this in time and ensuring
\begin{align}
\Bigg|  \frac{d}{dt} \int g^{ij}(u) \Big\langle \frac{\partial \psi}{\partial u^i}, \frac{\partial \psi}{\partial u^j} \Big\rangle_{\mathbb{R}^D} \sqrt{\text{det}(g)} \bigwedge_k du^k\Bigg| > \delta > 0
\end{align}
for suffciently large $\delta$. Under basic regularity conditions, we can exchange differentiation and integration. What follows is a very simple application of differentiating an inner product (in fact, this is not necessary but it makes automatic differentiation easier)
\begin{align}
& \Bigg|  \int_{\Sigma} \frac{d}{dt} g^{ij}(u) \Big\langle \frac{\partial \psi}{\partial u^i}, \frac{\partial \psi}{\partial u^j} \Big\rangle_{\mathbb{R}^D} \sqrt{\text{det}(g)} \bigwedge_k du^k\Bigg| 
\\
& = \Bigg|  \int_{\Sigma} g^{ij}(u) \Big( \Big\langle  \frac{\partial}{\partial t} \frac{\partial \psi}{\partial u^i}, \frac{\partial \psi}{\partial u^j} \Big\rangle_{\mathbb{R}^D}  + \Big\langle  \frac{\partial \psi}{\partial u^i}, \frac{\partial}{\partial t} \frac{\partial \psi}{\partial u^j} \Big\rangle_{\mathbb{R}^D}  \Big) \sqrt{\text{det}(g)} \bigwedge_k du^k\Bigg| 
\end{align}
Thus, we will set our geometric flow loss
\begin{align}
\label{eqn:harmonic_loss_app}
\mathcal{L} = \Bigg| \text{positive constant} - \Bigg|  \int_{\Sigma} g^{ij}(u) \Big( \Big\langle  \frac{\partial}{\partial t} \frac{\partial \psi}{\partial u^i}, \frac{\partial \psi}{\partial u^j} \Big\rangle_{\mathbb{R}^D}  + \Big\langle  \frac{\partial \psi}{\partial u^i}, \frac{\partial}{\partial t} \frac{\partial \psi}{\partial u^j} \Big\rangle_{\mathbb{R}^D}  \Big) \sqrt{\text{det}(g)} \bigwedge_k du^k\Bigg|  \Bigg|.
\end{align}

\vspace{2mm}

\textbf{Theorem 5 (slightly informal).} The loss function of \ref{eqn:harmonic_loss_app} has nontrivial dependence on Ricci curvature.

\vspace{2mm}

\textit{Proof.} We begin with the functional and apply integration by parts. Assuming compact support over $\Sigma$ such that boundary terms vanish, and applying Clairaut's theorem, we have
\begin{align}
    &\int_{\Sigma} g^{ij}(u) \Big( \Big\langle \frac{\partial}{\partial t} \frac{\partial \psi}{\partial u^i}, \frac{\partial \psi}{\partial u^j} \Big\rangle_{\mathbb{R}^D} + \Big\langle \frac{\partial \psi}{\partial u^i}, \frac{\partial}{\partial t} \frac{\partial \psi}{\partial u^j} \Big\rangle_{\mathbb{R}^D} \Big) \sqrt{\det g} \bigwedge_k du^k\nonumber \\
    & \propto \int_{\Sigma} g^{ij}(u) \Big\langle \frac{\partial}{\partial u^i} \Big(\frac{\partial \psi}{\partial t}\Big), \frac{\partial \psi}{\partial u^j} \Big\rangle_{\mathbb{R}^D} \sqrt{\det g} \bigwedge_k du^k\nonumber \\
    & \propto \int_{\Sigma} \langle \partial_t \psi, \Delta_g \psi \rangle \sqrt{\det g} \bigwedge_k du^k\nonumber \\
    & \propto  \int_{\Sigma} \langle \psi, \Delta_g (\partial_t \psi) \rangle \sqrt{\det g} \bigwedge_k du^k .
\end{align}
The above is defined element-wise, i.e.
\begin{align}
\Delta_g \psi = (\Delta_g \psi^1, \dots, \Delta_g \psi^D) .
\end{align}
We decompose the flow velocity field into its tangential and normal components \cite{doCarmo}, $\partial_t \psi = W + V^{\perp}$, where $W = W^\alpha \partial_\alpha \psi \in T\Sigma$ is tangential and $V^{\perp}$ is normal. The ambient Laplacian of the tangential component $W$ is
\begin{align}
    \Delta_g W = g^{ij} \overline{\nabla}_i \overline{\nabla}_j (W^\alpha \partial_\alpha \psi).
\end{align}
Applying the ambient connection $\overline{\nabla}$ and the Gauss formula $\overline{\nabla}_j \partial_\alpha \psi = {\Gamma_{j \alpha}}^\beta \partial_\beta \psi + \Pi_{j \alpha}$, the tangential projection onto $T\Sigma$ is governed by the intrinsic covariant derivative and the Weingarten map $A$,
\begin{align}
    (\Delta_g W)^{\top} = g^{ij} (\nabla_i \nabla_j W^\alpha) \partial_\alpha \psi - g^{ij} W^\alpha A_{\Pi_{i \alpha}}(\partial_j).
\end{align}
It can be noted $W$ is tangential but $\Delta_g W$ is not. We relate the connection Laplacian $\Delta_g W^\alpha = g^{ij} \nabla_i \nabla_j W^\alpha$ to the Hodge-de Rham Laplacian $\Delta_H$ using the Bochner-Weitzenb\"{o}ck identity \cite{doCarmo},
\begin{align}
    \Delta_H W^\alpha = -\Delta_g W^\alpha + R_{\beta}{}^\alpha{} W^\beta,
\end{align}
where $R_{\beta}{}^\alpha{}$ is the intrinsic Ricci curvature tensor. Rearranging this identity allows us to substitute the connection Laplacian,
\begin{align}
    g^{ij} \nabla_i \nabla_j W^\alpha = -\Delta_H W^\alpha + R_{\beta}{}^\alpha{} W^\beta.
\end{align}
Substituting this relation back into the tangential projection, we obtain
\begin{align}
    (\Delta_g W)^{\top} = \big( -\Delta_H W^\alpha + R_{\beta}{}^\alpha{} W^\beta \big) \partial_\alpha \psi - g^{ij} W^\alpha A_{\Pi_{i \alpha}}(\partial_j).
\end{align}
Reconstructing the right-hand side of the integral with this decomposition incorporates the Ricci curvature,
\begin{align}
    &\int_{\Sigma} \langle \psi, \Delta_g (\partial_t \psi) \rangle \sqrt{\det g} \bigwedge_k du^k\nonumber \\
    &\quad= \int_{\Sigma} \Big\langle \psi, \Delta_g (W + V^{\perp}) \Big\rangle \sqrt{\det g} \bigwedge_k du^k\nonumber \\
    &\quad= \int_{\Sigma} \bigg\langle \psi, \big( {-\Delta_H W^\alpha + R_{\beta}{}^\alpha{} W^\beta} \big) \partial_\alpha \psi - g^{ij} W^\alpha A_{\Pi_{i \alpha}}(\partial_j) + (\Delta_g V^{\perp})^{\top} + \text{other terms} \bigg\rangle \sqrt{\det g} \bigwedge_k du^k .
\end{align}
This confirms the nontrivial dependence on curvature. The proof is slightly informal because we do not show the Ricci curvature term is not canceled out.

$ \square $

\section{Additional experimental details}

Throughout our experiments, we will use the Riemannian metric of the sphere on occassion, which is
\begin{align}
ds^2 = R^2 \Bigg( du_1^2 +  \sin^2 u_1 \Bigg( du_2^2 + \sin^2 u_2 \Bigg( \hdots \Bigg( du_{n-1}^2  + du_n^2 \sin^2 u_{n-1} \Bigg) \Bigg) \Bigg)  \Bigg) .
\end{align}

\subsection{Navier-Stokes on the sphere}

In this subsection, we consider learning an encoder-decoder sequence on the Navier-Stokes equation
\begin{align}
\begin{cases}
\partial_t v(x,t) + (v(x,t) \cdot \nabla) v(x,t) - \nu \nabla^2 v(x,t) = \frac{-1}{\rho} \nabla p(x,t) 
\\
\nabla \cdot v = 0
\\
(x,t) \in \mathcal{X} \times [0,T] \subseteq \mathbb{R}^2 \times [0,T]
\end{cases} .
\end{align}
Specifically, we learn the vorticity over an $80 \times 80$ mesh in $\mathbb{R}^2$. Data was created using the Galerkin/Fourier spectral solver of \cite{crewsincompressibleNS} 
\cite{crewsnavierstokes}. We will consider the special case of Ricci flow on the expanding sphere in this setting, and our manifold takes the form
\begin{align}
\mathcal{E}_{\mathbb{S}^d} = \frac{ (u + \xi) \sqrt{50^2 + 50(d-1)t}}{ || u + \xi ||_2 } .
\end{align}
Here, $\xi \sim 0.375 \cdot \mathcal{N}^d(0,1)$ is a vector of Gaussian noise. By normalizing the point in space after the noise, we ensure the noise lies within the manifold and not around it. We have chosen $d=100$. Here, we have denoted $d$ the extrinsic dimension, since this manifold is not constructed with a parametric map. Our neural networks in this experiment are $u = u_{\theta_u}$, $\mathcal{D} = \mathcal{D}_{\theta_{\mathcal{D}}}$, and $\mathcal{S} = \mathcal{S}_{\theta_{\mathcal{S}}}$. In the test setting, we take
\begin{align}
\mathcal{E}_{\mathbb{S}^d} = \frac{( u + \EX_{\xi} [\xi] ) \sqrt{50^2 + 50(d-1)t}}{  || u + \EX_{\xi} [ \xi ] ||_2  } = \frac{ u \sqrt{50^2 + 50(d-1)t}}{ || u  ||_2  } .
\end{align}
We take a convolutional neural network for $u$ and a vanilla multilayer perceptron for decoder $\mathcal{D}$, which we found outperformed a transpose convolutional (deconvolutional) network. Our decoder used severe dropout (3 iterations with $p=0.325$ and 1 iteration with $p=0.225$). We chose a fairly high learning rate $5\text{e}{-4}$ and fixed it. 

\vspace{2mm}

For our baselines, we only construct $u,\mathcal{D}$. We construct an autoencoder with an unrestricted latent stage. For two scenarios, we add: (1) no noise in the first setting; (2) the same exact noise we introduced into our manifold in the unrestricted latent space.

\subsection{Burger's equation}
\label{sec:burgers_experiments}

\textbf{\begin{figure}[!h]
  \vspace{0mm}
  \centering
  \includegraphics[scale=0.65]{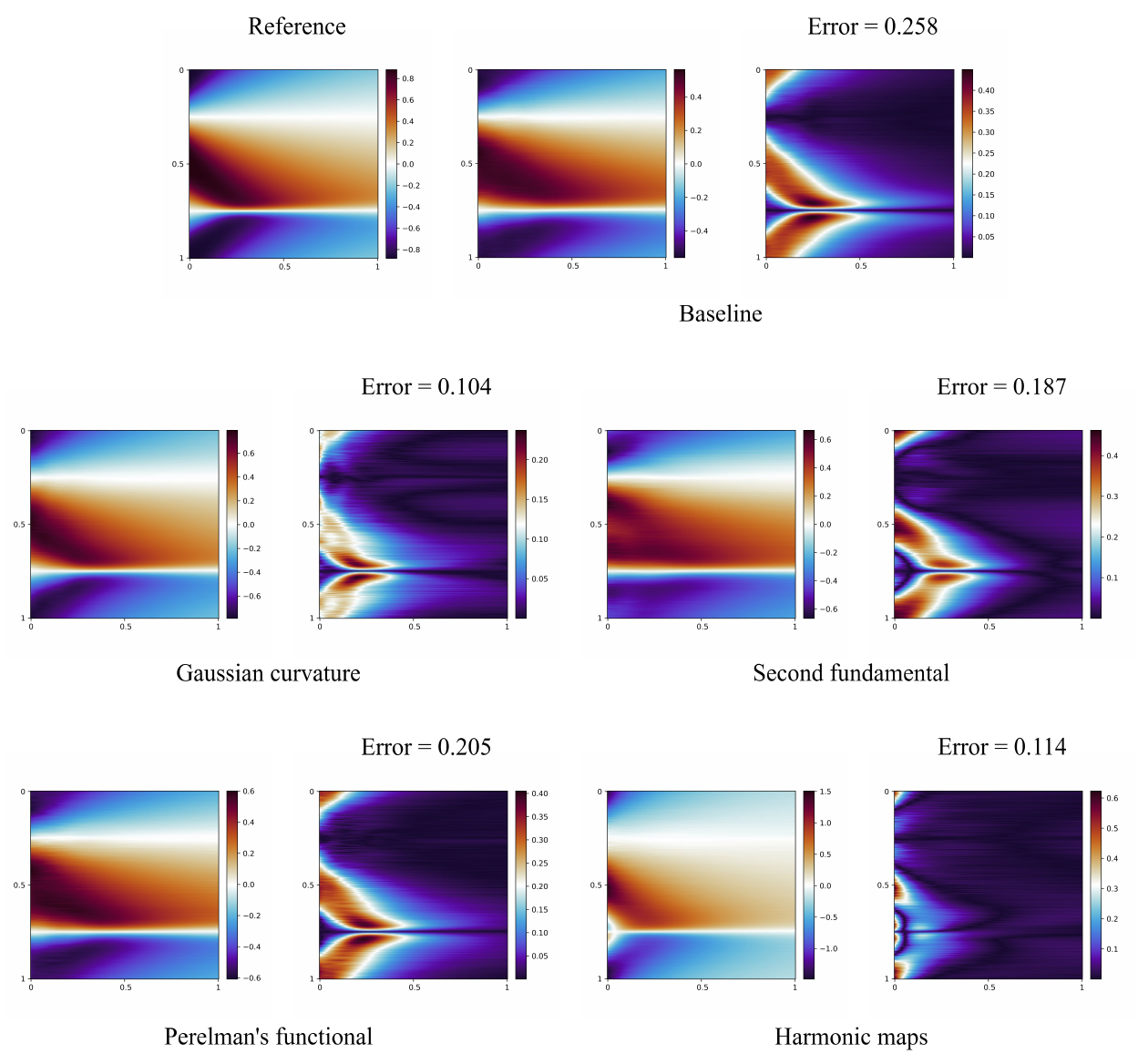}
  \caption{We illustrate our results on a single example of Burger's equation in the variational setting. The OOD scenario here is noise with $\sigma=1.5$ is injected into 101 locations of the initial condition. We highlight this is only on a single example, so refer to Table \ref{tab:method_errors} for results across numerous data. We also highlight our Gaussian curvature experiment has intrinsic dimension $d=2$, and extrinsic dimension $D=3$. All other experiments have intrinsic dimension $d=3$, and extrinsic dimension $D=5$. All hyperparameters are held as constant as possible across training, although the hyperparameter that is given a little bit of freedom is the loss scaling coefficients since each loss function has different terms.}
\label{fig:geoflow_burgers}
\end{figure}}

The remainder of our experiments will be restricted to one dataset. This is to ensure we demonstrate all facets of our methods across a consistent datatset for comparison, and this narrow view allows us to deep-dive into the outcomes of our methods.

\vspace{2mm}

We will consider the Burger's equation
\begin{gather}
\partial_t \phi = \nu \partial_{xx} \phi  - \phi \partial_x \phi
\\
\text{such that}
\\
(x,t) \in \mathcal{X} \times [0,T] = [0,1]\times [0,1]
\\
\phi(x,0) = \sum_{i=1}^3 \alpha_i \cos^{2i-1}(2 \pi x) , \alpha_i \sim U[-1,1].
\end{gather}
For our baseline, we will consider extended encoder methods, meaning we take
\begin{align}
\tilde{\phi}_t \approx \mathcal{D} \circ \mathcal{E} \circ u (\tilde{\phi}_0, t)
\end{align}
but without any restrictions on $\mathcal{E}$. We consider this because it is a more similar architecture to ours, so we will investigate the geometric flow effects only. Also, this architecture generally has better out-of-distribution robustness than a baseline encoder $\mathcal{D} \circ \mathcal{E}$ (see \cite{gracyk2025variationalautoencoderslatenthighdimensional}).

\vspace{2mm}

We will exclusively consider vanilla multilayer perceptron architectures for this experiment, since this is the most typical architecture. Note that the architecture of \cite{wang2023expertsguidetrainingphysicsinformed} works well as well, but this architecture is a bit specialized, so we will not examine this. We will take all neural networks to be of the form
\begin{align}
\text{MLP}(x) = W^{L} ( \sigma \circ (  W^{L-1} \sigma \circ \hdots \circ (\sigma \circ (W^1 x + b^1 ) + b^2) + \hdots + b^{L-1} ) + b^L .
\end{align}
We will take
\begin{align}
\sigma(\cdot) = \text{tanh}(\cdot)
\end{align}
for all experiments in any network in which the PDE data is passed. We remark it is common to normalize layers sometimes, such as in \cite{banerjee2022restrictedstrongconvexitydeep}, but we will not do this. When we construct a metric $g$ with a neural network $g_{\theta_g}$, we will always take
\begin{align}
g(u,t) = (\mathcal{G}(u,t))^T \mathcal{G}(u,t)
\end{align}
to ensure positive semi-definiteness.

\vspace{2mm}

In table \ref{tab:method_errors}, we remark this is quite a large sample size. Data is gathered so that we examine the aggregate relative $L^1$ error every 10-th discrete time point in the 201-point mesh in time. All of our methods are trained until roughly $2\text{e}{-4}$ decoded loss (MSE on the PDE learning) or until convergence (we remark the Gaussian curvature experiment was not trained as far, which was favorable and is likely attributable to its dimension differences, but in general we did not find this extremely important). We choose KL regularization coefficient $\beta=0.001$. We remark that high batch size and low learning rate ($2\text{e}{-4}$) led to consistent results among retraining, especially among the baseline method (see Figure \ref{fig:errors_consistent}), thus it is a meaningful benchmark. The seed is set such that all data is the same among the comparisons, both the PDE data and the OOD data.

\vspace{2mm}

\textbf{Gaussian curvature experiment.} We list details for our Gaussian curvature experiment. Empirically, we tried the change of variables strategy as in Appendix \ref{app:gauss_curv_sec}. This loss function is extremely nontrivial to optimize, and this method is extremely suboptimal. Instead, we found our division by a nonzero function method much simpler and more effective.

\vspace{2mm}

\textbf{Second fundamental form with its proxies.} We found good empirical results with the addition of a further diagonal term based on $\Pi$ and a small adjustment of dimensional constants. We preserve the core structure of the theoretical setup, thus our results are consistent with what we established in Appendix \ref{app:second_FF_sec}. Also, we have approximated our Hessian using a Jacobian transpose Jacobian. This is a known approximation, and it reduces computational expense quite significantly. Here, we have initialized the metric at $t=0$ with $g_{\text{sphere}}$.

\vspace{2mm}

\textbf{Perelman's functional experiment.} For simplicity, we will set our baseline metric as the sphere $g_0 = g_{\text{sphere}}$. The sphere has scalar curvature
\begin{align}
R = \frac{d(d-1)}{r^2} .
\end{align}
Here, $r$ is the radius.

\vspace{2mm}

\textbf{Harmonic map experiment.} For our initial manifold $\mathcal{M}$, we will map our data to any point in Euclidean space $\mathbb{R}^d$ and we will normalize this data such that
\begin{align}
u \rightarrow \frac{r \times u}{||u||_2}
\end{align}
for some radius $r$. $\psi$ is our parametric map such that
\begin{align}
\psi_{\theta_{\psi}} :  \frac{r \times u}{||u||_2} \rightarrow \mathcal{N}_t \subset \mathbb{R}^D .
\end{align}
We set $r=2$. We use a geometric flow loss of the form
\begin{align}
\Bigg| 0.5 -   \int_{\Sigma} \frac{d}{dt} \Big( g^{ij}(u) \Big\langle \frac{\partial \psi}{\partial u^i}, \frac{\partial \psi}{\partial u^j} \Big\rangle_{\mathbb{R}^D}  \Big) \sqrt{\text{det}(g)} \bigwedge_k du^k\Bigg|^2 + \EX \Bigg[ \Bigg| \ \Bigg| \Bigg| \psi(u,0) \Bigg| \Bigg|_2 - 5 \ \Bigg|^2 \Bigg] ,
\end{align}
which yielded both consistent and favorable results.

\section{Additional figures}

\newpage

\begin{figure}[!h]
  \vspace{0mm}
  \centering
  \includegraphics[scale=0.78]{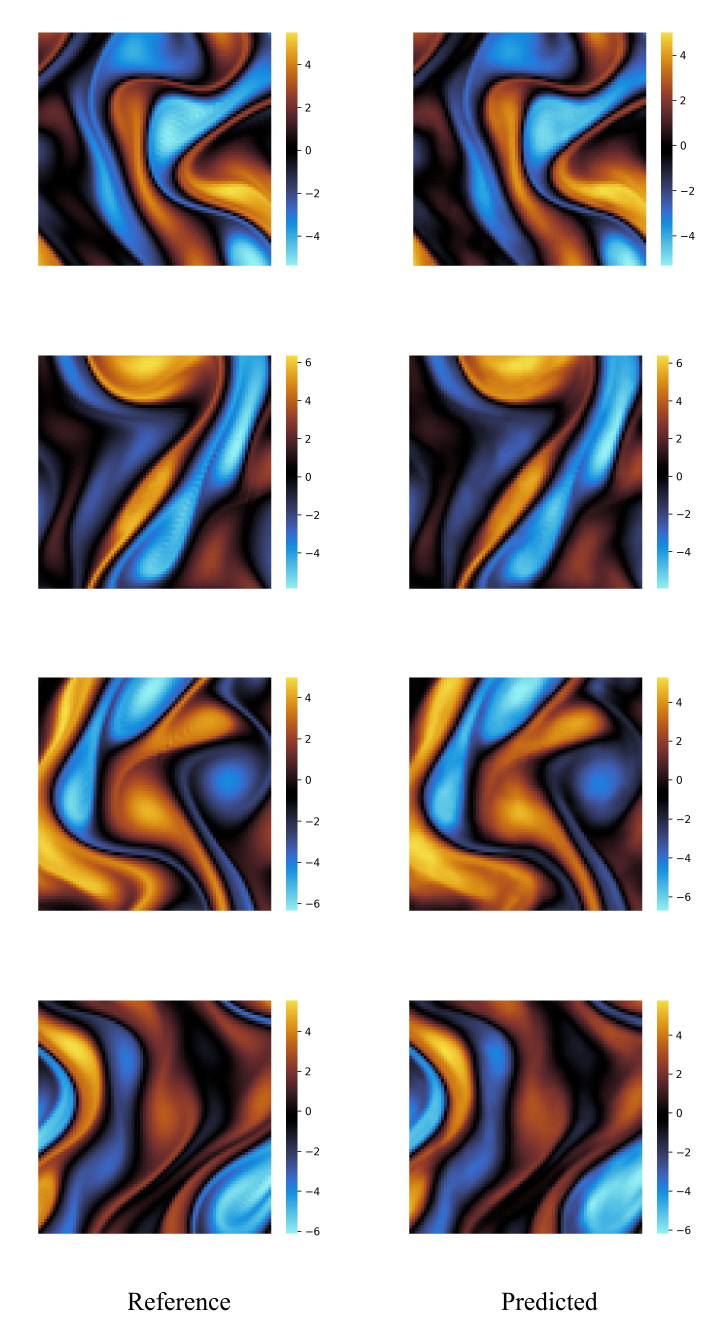}
  \caption{We present uncurated (except the first) results on our heavily-regularized autoencoder on Navier-Stokes data using Ricci flow on the sphere. This figure illustrates in-distribution data with heavy dropout and latent Gaussian noise in the test setting, not training. We find a manifold with very highly-expanding volume element works best. We also remark training data heavily favored final time $t=T$ in amount of data, which this data reflects. We remark the baseline autoencoder actually performs similarly in this test stage (which is partially because the data can get learned better in the sense of lower training loss without the regularization), but it is in robustness inference where our method shines, as is demonstrated in Figure \ref{fig:navier_stokes_ood}. We also remark the regularization techniques with our method were found to be necessary, but not in a vanilla autoencoder.}
\label{fig:navier_stokes_indist}
\end{figure}

\textbf{\begin{figure}[htbp]
  \vspace{0mm}
  \centering
  \includegraphics[scale=0.75]{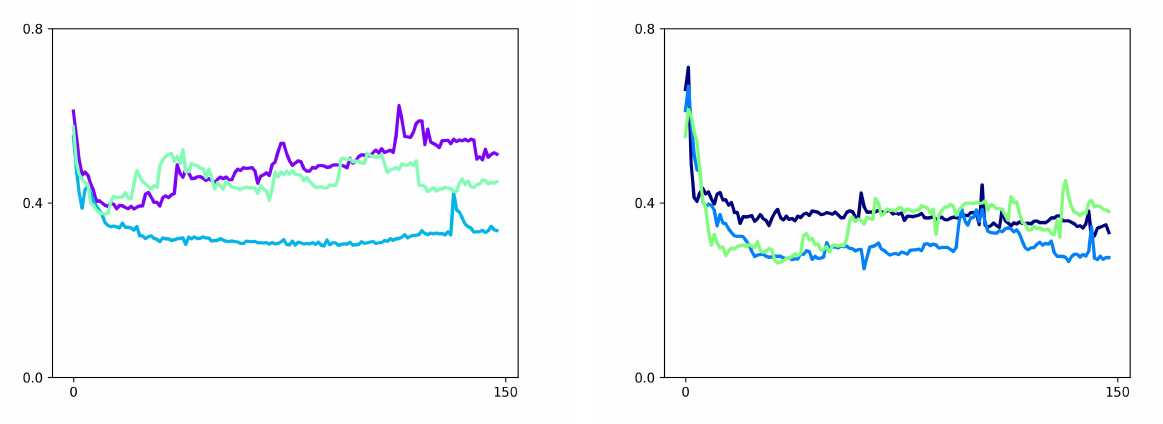}
  \caption{We compare a geometric flow of the form $\partial_t g = - \Lambda$ to the steady geometric flow of \cite{gracyk2025variationalautoencoderslatenthighdimensional}, which produces both a more canonical geometry and one that is of large measure. We plot 3 trajectories with severely noised initial condition and sample discrete relative $L^2$ error at 5 select times along the time interval. As we can see, the more optimal geometry has more favorable outcomes. }
\label{fig:geo_flow_comparison}
\end{figure}}

\textbf{\begin{figure}[htbp]
  \vspace{0mm}
  \centering
  \includegraphics[scale=0.75]{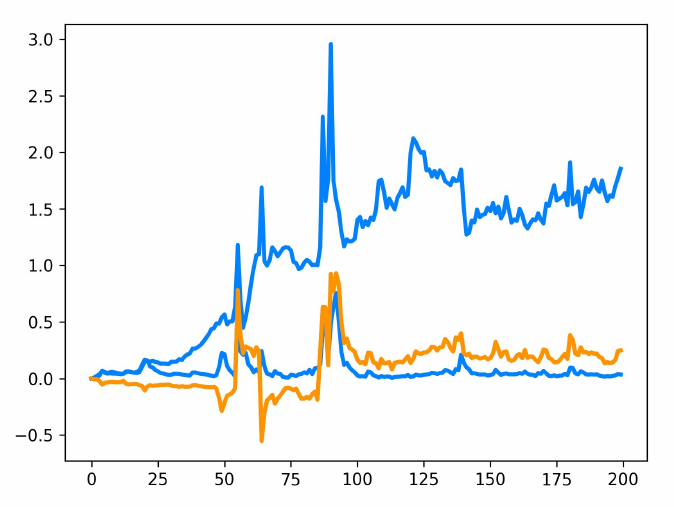}
  \caption{We provide empirical evidence that curvature-based flows help support nondegeneracy. This figure represents the means (over the batch, so with respect to both PDE data and time) of the $2 \times 2$ Riemannian metric learned in our Gaussian curvature experiment with respect to training iteration. As we can see, our metrics are over threefold larger than that of the \ref{fig:small_metric_vae} experiment at the same training iteration, which is mod 20. Blue is the diagonal, and orange is the off-diagonal.}
\label{fig:metric_gausscurv}
\end{figure}}

\textbf{\begin{figure}[htbp]
  \vspace{0mm}
  \centering
  \includegraphics[scale=0.74]{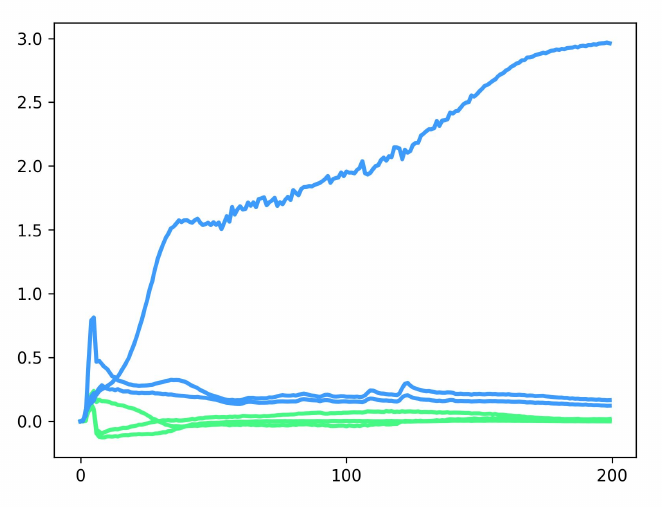}
  \caption{We provide a similar figure to that of \ref{fig:metric_gausscurv} but on our custom geometric flow based on the Riemannian decomposition with the inexpensive second fundamental form. It is noteworthy that we have initialized our metric as $g(u,0) = 3 g_{\text{unit sphere}}$ here, which helped slightly. Blue is the diagonal, which corresponds to $\Gamma_{ijkl}$, and green is the off-diagonal.}
\label{fig:metric_secondfundamental}
\end{figure}}

\textbf{\begin{figure}[htbp]
  \vspace{0mm}
  \centering
  \includegraphics[scale=0.74]{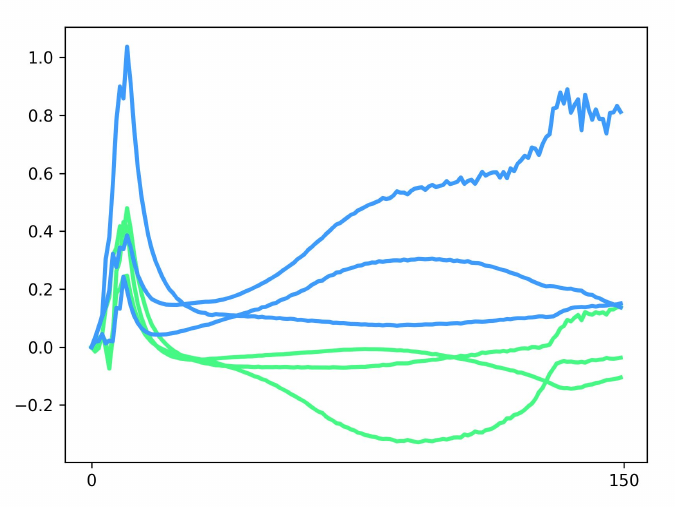}
  \caption{This figure goes hand-in-hand with \ref{fig:metric_secondfundamental} and metric belonging to the custom Riemannian decomposition experiment, but we instead initialize our metric at $t=0$ with $\cos(u^1) du^1 \otimes du^3$ in the $(3,1),(1,3)$ off-diagonal. Note that this helps support a nondegenerate off-diagonal, but the diagonal is not as large.}
\label{fig:metric_secondfundamental_nosphere}
\end{figure}}

\textbf{\begin{figure}[htbp]
  \vspace{0mm}
  \centering
  \includegraphics[scale=0.75]{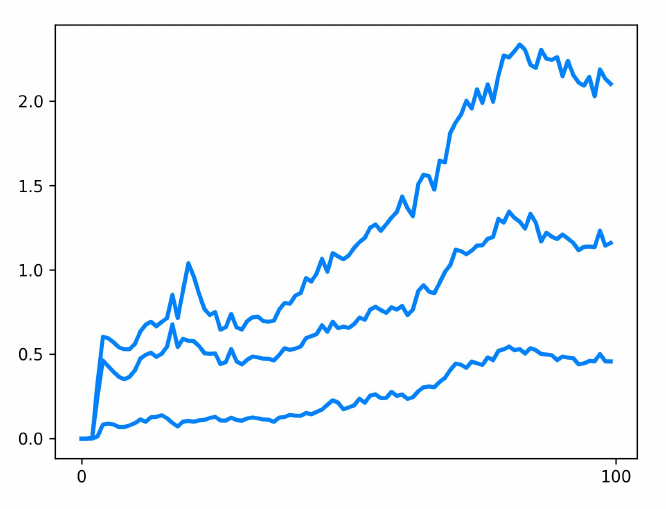}
  \caption{We provide a similar figure to that of \ref{fig:metric_gausscurv} but of our Perelman functional experiment with a conformally changed metric. As we can see, nondegeneracy is also supported in this experiment. Here, we have clipped $\psi$ (not its derivatives) to be sufficiently large, which we have found works well. Since our baseline metric here is diagonal (the sphere), we only plot the diagonal.}
\label{fig:metric_perelman}
\end{figure}}

\textbf{\begin{figure}[htbp]
  \vspace{0mm}
  \centering
  \includegraphics[scale=0.75]{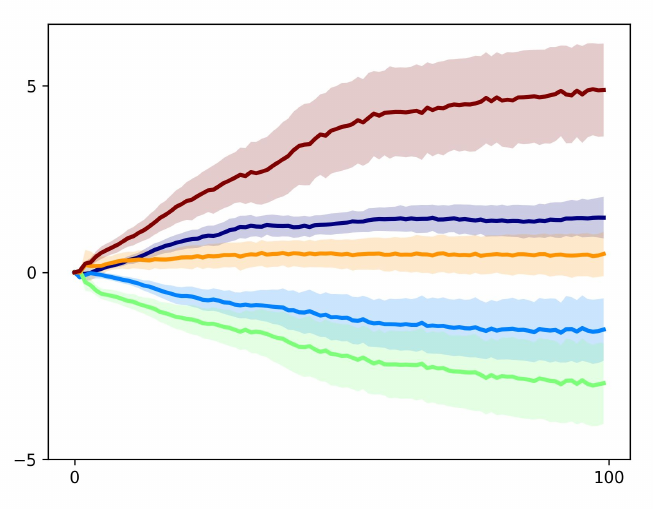}
  \caption{We plot the means and standard deviations of the parametric map in the harmonic map experiment with respect to training iteration mod 20. The standard deviation especially helps show nontriviality, since the manifold can be immersed anywhere in space.}
\label{fig:parametric_map_harmonic}
\end{figure}}

\textbf{\begin{figure}[!h]
  \vspace{0mm}
  \centering
  \includegraphics[scale=0.75]{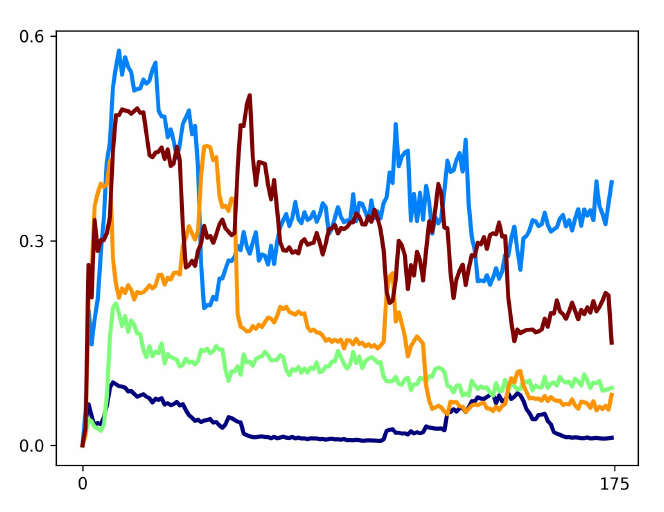}
  \caption{We plot the $5 \times 5$ mean over the batch (which includes over the time interval) diagonal elements of a learned Riemannian metric with a geometric flow of the form $\partial_t g(u,t) = -\Lambda(u,t)$ (both are positive semi-definite) in the variational setting on Burger's equation data with respect to training iteration mod 20. We do not need to plot the off-diagonal, as these scale according to the diagonal by Cauchy-Schwarz. As we can see, values of the metric are fairly small, and so the manifold is not of large measure. Decoded coefficient here was 100, with metric and flow coefficients of 1.}
\label{fig:small_metric_vae}
\end{figure}}

\textbf{\begin{figure}[!h]
  \vspace{0mm}
  \centering
  \includegraphics[scale=0.75]{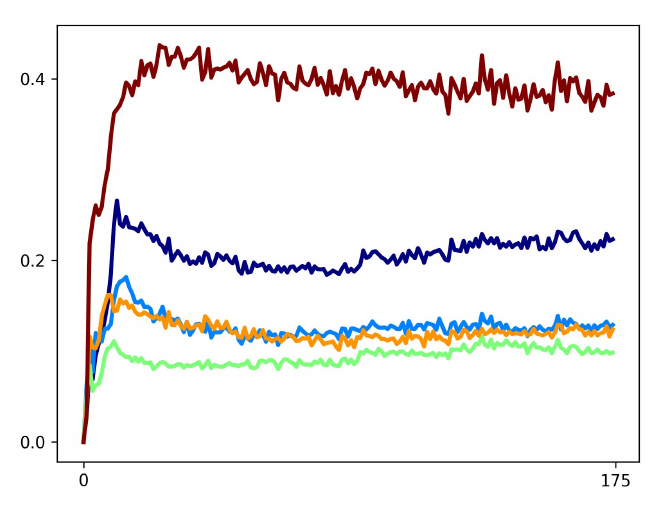}
  \caption{We plot the $5 \times 5$ mean over the batch (which includes over the time interval) diagonal elements of a learned Riemannian metric with a geometric flow of the form $\partial_t g(u,t) = -\Lambda(u,t)$ in the encoder-decoder setting on Burger's equation data with respect to training iteration mod 20. Decoded coefficient here was 10, with metric and flow coefficients of 1. }
\label{fig:small_metric_ae}
\end{figure}}

\begin{figure}[!h]
  \vspace{0mm}
  \centering
  \includegraphics[scale=0.58]{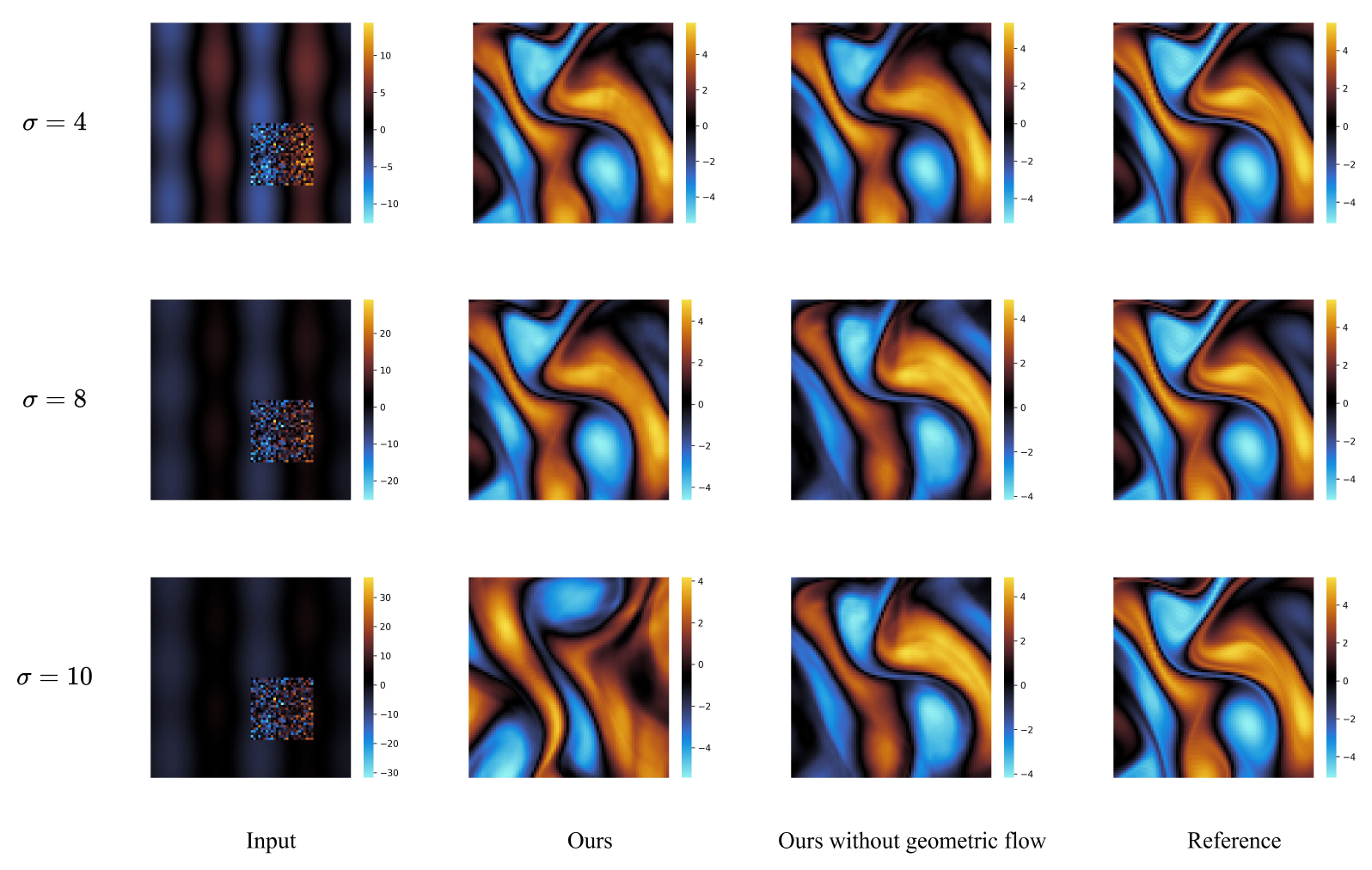}
  \caption{We demonstrate more uncurated (curated for complexity but not for inference quality) results on our Navier-Stokes experiment. We have found our method behaves significantly better with the geometric flow than without under reasonable noise. Eventually, if the OOD effect is so severe, our method outputs an entirely distinct solution to that that is the true data. We remark the OOD effect needs to be extremely severe to achieve this breaking point.}
\label{fig:navier_stokes_breakingpoint}
\end{figure}

\end{document}